\documentclass[preprint,3p,times]{elsarticle}
\usepackage{amssymb}

\usepackage{amsmath}

\usepackage{hyperref} 
\usepackage{marginnote}
\usepackage{xcolor}
\usepackage{amsthm}
\usepackage{float}
\usepackage{tikz}
\usepackage{booktabs}
\usepackage{pgfplots}
\pgfplotsset{compat = 1.3}
\usepackage{mathtools}
\usepackage{caption}
\usepackage{subcaption}

\biboptions{numbers,sort&compress}

\newtheorem{theorem}{Theorem}

\newtheorem{remark}{Remark}
\newtheorem{example}{Example}

\newcommand{\Ham}{\mathcal{H}}
\newcommand{\T}{\mathcal{T}}
\newcommand{\V}{\mathcal{V}}
\newcommand{\bfA}{\boldsymbol{A}}
\newcommand{\bfB}{\boldsymbol{B}}
\newcommand{\bfAhat}{\boldsymbol{\hat{A}}}
\newcommand{\bfBhat}{\boldsymbol{\hat{B}}}
\newcommand{\A}{\mathcal{A}}
\newcommand{\B}{\mathcal{B}}
\newcommand{\C}{\mathcal{C}}
\newcommand{\D}{\mathcal{D}}
\newcommand{\p}{\boldsymbol{p}}
\newcommand{\q}{\boldsymbol{q}}
\newcommand{\bfP}{\boldsymbol{P}}
\newcommand{\bfQ}{\boldsymbol{Q}}

\newcommand{\0}{\boldsymbol{0}}

\newcommand{\G}{\mathcal{G}}
\newcommand{\frakg}{\mathfrak{g}}
\newcommand{\bfT}{\boldsymbol{T}}
\newcommand{\e}{\boldsymbol{e}}
\newcommand{\bomega}{\boldsymbol{\omega}}
\newcommand{\diff}{\mathrm{d}}
\newcommand{\bfM}{\boldsymbol{M}}

\definecolor{fk4blue}{rgb}{0.00000,0.46275,0.68627}
\journal{Computer Physics Communications}

\begin{document}

\begin{frontmatter}

%% Title, authors and addresses

%% use the tnoteref command within \title for footnotes;
%% use the tnotetext command for theassociated footnote;
%% use the fnref command within \author or \address for footnotes;
%% use the fntext command for theassociated footnote;
%% use the corref command within \author for corresponding author footnotes;
%% use the cortext command for theassociated footnote;
%% use the ead command for the email address,
%% and the form \ead[url] for the home page:
%% \title{Title\tnoteref{label1}}
%% \tnotetext[label1]{}
%% \author{Name\corref{cor1}\fnref{label2}}
%% \ead{email address}
%% \ead[url]{home page}
%% \fntext[label2]{}
%% \cortext[cor1]{}
%% \affiliation{organization={},
%%             addressline={},
%%             city={},
%%             postcode={},
%%             state={},
%%             country={}}
%% \fntext[label3]{}

\title{Hessian-free force-gradient integrators}

\author[1]{Kevin Schäfers\corref{cor1}}
\ead{schaefers@math.uni-wuppertal.de}
\author[2,3]{Jacob Finkenrath}
\ead{j.finkenrath@cern.ch}
\author[1]{Michael Günther}
\ead{guenther@math.uni-wuppertal.de}
\author[2]{Francesco Knechtli}
\ead{knechtli@physik.uni-wuppertal.de}

\cortext[cor1]{Corresponding author.}

\affiliation[1]{organization={Institute of Mathematical Modelling, Analysis and Computational Mathematics (IMACM), Chair of Applied and Computational Mathematics, Bergische Universität Wuppertal},%Department and Organization
            addressline={Gaußstraße 20}, 
            city={Wuppertal},
            postcode={42119}, 
            %state={},
            country={Germany}}
\affiliation[2]{organization={Department of Physics, Bergische Universität Wuppertal},%Department and Organization
            addressline={Gaußstraße 20}, 
            city={Wuppertal},
            postcode={42119}, 
            %state={},
            country={Germany}}
\affiliation[3]{organization={Department of Theoretical Physics},%Department and Organization
            addressline={European Organization for Nuclear
Research, CERN}, 
            city={Geneve},
            postcode={1211}, 
            %state={},
            country={Switzerland}}

\begin{abstract}
We propose a new framework of Hessian-free force-gradient integrators that do not require the analytical expression of the force-gradient term based on the Hessian of the potential. Due to that the new class of decomposition algorithms for separable Hamiltonian systems with quadratic kinetic energy may be particularly useful when applied to Hamiltonian systems where an evaluation of the Hessian is significantly more expensive than an evaluation of its gradient, e.g.~in molecular dynamics simulations of classical systems. Numerical experiments of an N-body problem, as well as applications to the molecular dynamics step in the Hybrid Monte Carlo (HMC) algorithm for lattice simulations of the Schwinger model and Quantum Chromodynamics (QCD) verify these expectations.
\end{abstract}

\begin{keyword}
%% keywords here, in the form: keyword \sep keyword
Geometric integration \sep Hamiltonian systems \sep Splitting methods \sep Order conditions \sep Hybrid Monte Carlo \sep Lattice Quantum Chromodynamics

%% PACS codes here, in the form: \PACS code \sep code

%% MSC codes here, in the form: \MSC code \sep code
%% or \MSC[2008] code \sep code (2000 is the default)
\MSC[2020]
81V05 % Strong interaction, including quantum chromodynamics 
\sep 
65P10 % Numerical methods for Hamiltonian systems including symplectic integrators
\sep
65L05 % Numerical methods for initial value problems involving ordinary differential equations
\sep
65L20 % Stability and convergence of numerical methods for ordinary differential equations
37N20 % Dynamical systems in other branches of physics (quantum mechanics, general relativity, laser physics)

\end{keyword}
\end{frontmatter}

%% \linenumbers

%% main text
\section{Introduction}
The numerical integration of Hamiltonian systems imposes challenging demands on the numerical integration scheme. Particularly, the phase space $(\p,\q)$ with generalized coordinates $\q$ and conjugate momenta $\p$ is a symplectic manifold, giving rise to the need of \emph{geometric numerical integration} \cite{HairerLubichWanner}. 

In this paper, we consider separable Hamiltonian systems $\Ham(\p,\q) = \T(\p) + \V(\q)$ with kinetic energy of the form 
\begin{align}\label{eq:T}
\T(\p) = \tfrac{1}{2}\p^\top \bfM^{-1} \p,
\end{align}
where $\bfM$ is a constant symmetric positive definite matrix.
This particular structure has a wide range of applications in the field of many-body problems, e.g.~classical mechanics \cite{goldstein2002classical}, quantum mechanics \cite{griffiths2018introduction}, statistical mechanics \cite{gibbs1902elementary} and lattice field theories \cite{duane1987hybrid}. In the latter case, the phase space is the cotangent bundle over a base space that is a Lie group manifold. Moreover, any system of second order ODEs $\ddot{y} = f(y)$ is reduced to the equations of motion under consideration in this work.

There exist many approaches for geometric numerical integration of Hamiltonian systems. Recently, symplectic generalized additive Runge--Kutta (GARK) schemes \cite{gunther2023symplectic,schafers2023symplectic} have been developed. An alternative approach based on the variational principle has been discussed in \cite{zanna2020discrete}. Both approaches are restricted to Euclidean space.
An extension of symplectic partitioned Runge--Kutta methods to the Lie group setting has been presented in \cite{bogfjellmo2016high}. However, the general framework does not include time-reversibility, another important feature of Hamiltonian systems. 
For separable Hamiltonian systems, splitting methods allow for the derivation of explicit geometric numerical integration schemes of arbitrarily high convergence via composition techniques \cite{suzuki1990fractal,yoshida1990construction,omelyan2002construction}. 

By exploiting the special structure of the kinetic energy $\mathcal{\T}(\p)$, \emph{force-gradient integrators} (FGIs) \cite{omelyan2003symplectic,kennedy2009force} allow for an efficient computational process. FGIs can be regarded as splitting methods, applied to a modified potential, resulting in more accurate numerical approximations with respect to the original system.
One drawback of FGIs is the necessity of deriving the so-called \emph{force-gradient term} (FG-term), containing the Hessian of the potential $\V(\q)$.
To our knowledge, the idea of using the Hessian of the potential $\V(\q)$ to enhance splitting methods dates back to \cite{Takahashi1984,ROWLANDS1991235} and has been further investigated in \cite{Lopez_Hessian1997}. Recently, force-gradient integrators have been successfully developed in the more general context of port-Hamiltonian (pH) systems to derive higher order schemes fulfilling the dissipation inequality of pH systems and thus breaking the order-two limit for splitting methods~\cite{mönch2024commutatorbased}.
Another possible drawback is the evaluation cost of the FG-term. In molecular dynamics simulations of classical systems, for example, the evaluation typically is 2-3 times more expensive than a usual force evaluation \cite{omelyan2002construction}. One can overcome this issue by approximating the FG-term, as proposed in \cite{wisdom1996symplectic} and further investigated in \cite{Hairer_McLachlan_Skeel_2009} for a particular FGI. In the context of lattice quantum chromodynamics (QCD), the approximation of the FG-term was initially applied to another FGI in \cite{yin2011improving}. The adapted FGI from \cite{yin2011improving} has been applied to the two-dimensional Schwinger model where we implemented the FG-term analytically \cite{shcherbakov2017adapted} and in large scale lattice QCD simulations, see e.g.~\cite{Jung:2017xef,Finkenrath:2022eon,Finkenrath:2023sjg}. Numerical results highlight that the Hessian-free variant results in a more efficient computational process. 

In this paper we will generalize the idea of approximating the FG-term for the entire class of FGIs introduced in \cite{omelyan2003symplectic}. This will be applied not only in case of lattice field theories, but in general for all separable Hamiltonian systems with kinetic energy of the form \eqref{eq:T}, where the new class of \emph{Hessian-free force-gradient integrators} can be utilized. We will discuss the geometric properties, a refined error analysis of the approximation, resulting in explicit formulae for the error terms, as well as a backward error analysis giving insight into the long-time behavior of Hessian-free FGIs.

The paper is organized as follows. In Section \ref{sec:Hamiltonian_mechanics}, we will briefly introduce Hamiltonian mechanics and its demands on the numerical integration scheme. Section \ref{sec:Hessian-free_FGI} introduces the new class of Hessian-free force-gradient integrators. 
In Section \ref{sec:integrator_derivation}, Hessian-free force-gradient integrators with up to eleven stages are derived. Here, (weighted) norms of the leading error coefficients are defined as aggregated functions of the principal error term. By performing a global minimization of this aggregated function, we obtain optimal sets for the integrator coefficients.
In Section \ref{sec:Numerical_Results}, numerical results for three different test examples are discussed: a) the outer solar system, a $N$-body problem with $\bfM \neq \mathrm{Id}$, b) the two-dimensional Schwinger model where we can compare the performance of exact FGIs and Hessian-free FGIs as the analytical expression of the FG-term is available \cite{shcherbakov2017adapted}, and c) four-dimensional gauge field simulations in lattice QCD  with two heavy Wilson fermions. Here, the physical degrees of freedoms are elements of the non-Abelian matrix Lie group $\mathrm{SU}(3)$ and all tests are performed using the openQCD implementation \cite{openQCD,LUSCHER2013519}. The paper concludes with a summary and outlook for future research.

\section{Hamiltonian mechanics and geometric numerical integration}\label{sec:Hamiltonian_mechanics}
In this section, we will briefly introduce Hamiltonian mechanics on matrix Lie groups. Moreover we will state the demands on the numerical integration scheme for separable Hamiltonian systems 
\begin{equation}\label{eq:Hamiltonian_general}
    \Ham(\p,\q) = \T(\p) + \V(\q) = \tfrac{1}{2}\p^\top \bfM^{-1}\p + \V(\q)
\end{equation}
with constant symmetric positive definite matrix $\bfM$. Without loss of generality \footnote{with the Cholesky decomposition $\bfM \coloneqq \boldsymbol{L}^\top \boldsymbol{L}$ and transformed momenta $\tilde{\boldsymbol{p}} \coloneqq \boldsymbol{L}\p$, any system \eqref{eq:Hamiltonian_general} can be written as a system of the form \eqref{eq:Hamiltonian}}, we consider a separable Hamiltonian system of the form 
\begin{equation}\label{eq:Hamiltonian}
    \Ham(\p,\q) = \T(\p) + \V(\q) = \tfrac{1}{2} \langle \p,\p \rangle + \V(\q).
\end{equation}

\subsection{Hamiltonian mechanics}
We consider a phase space $(\p,\q) \in T^* \G$ where $T^* \G$ denotes the cotangent bundle over a base space that is a $d$-dimensional matrix Lie group manifold $\G$ and whose fibers are isomorphic to its Lie algebra $\frakg$. The linear space $\frakg$ has a basis, consisting of \emph{generators} $\bfT_i$, $i=1,\ldots,d$. For matrix Lie groups, there exists a matrix representation $\bfQ$ of the Lie group element $\q \in \G$ and $\bfP:= p^i \bfT_i \in \frakg$ of the momentum $\p$. The generators $\bfT_i$ are linked to the Lie group elements $\bfQ$ via the right-invariant linear differential operator $\e_i$ whose action on $\bfQ$ is defined by 
\begin{equation}\label{eq:diff_operator_e}
    \e_i(\bfQ) = - \bfT_i \bfQ.
\end{equation}
The operator \eqref{eq:diff_operator_e} can be regarded as a generalization of the vector fields $\partial/\partial\q$ in the Lie group space. The cotangent space has a natural symplectic structure $\bomega = -\diff\p$ that is closed, $\diff\bomega = 0$. The fundamental two-form $\bomega$ defines the Poisson bracket of two arbitrary Hamiltonian vector fields $\bfAhat,\bfBhat$ corresponding to zero-forms $\bfA,\bfB$ as $\{\bfA,\bfB\} = -\bomega(\bfAhat,\bfBhat)$. Hamilton's equations are most naturally expressed in terms of the Lie derivative operators 
\begin{equation*}
    \hat{\V} := \{\V, \cdot \} \quad \text{and} \quad \hat{\T} := \{\T,\cdot\},
\end{equation*}
resulting in the additively partitioned system of ordinary differential equations (ODEs)
\begin{align*}
    \begin{pmatrix}
        \dot{\bfP} \\ \dot{\bfQ}
    \end{pmatrix} &= \begin{pmatrix} \{\V, \bfP\} \\ \{ \V,\bfQ \} \end{pmatrix} + \begin{pmatrix}
        \{\T,\bfP\} \\ \{\T,\bfQ\}
    \end{pmatrix} = \begin{pmatrix} \{\V,\bfP\} \\ \0\end{pmatrix} + \begin{pmatrix} \0 \\ \{\T,\bfQ\} \end{pmatrix} = \begin{pmatrix} \hat{\V} \bfP \\ \0 \end{pmatrix} + \begin{pmatrix} \0 \\ \hat{\T}\bfQ \end{pmatrix}.
\end{align*}
For matrix Lie groups, there exist structure constants $c_{jk}^i$ satisfying $[\bfT_j,\bfT_k] = c_{jk}^i \bfT_i$. They occur in the Lie derivative operator
\begin{equation}\label{eq:kinetic_vf}
    \hat{\T} = p^i \e_i + c_{ki}^j p_j p^k \frac{\partial}{\partial p_i}
\end{equation}
of the kinetic part. For Abelian Lie groups (e.g. $\mathrm{U}(1)$ and $\mathrm{SO}(2)$), it holds $c_{jk}^i = 0$ for all $i,j,k$. In case of semisimple\footnote{a semisimple Lie group is a non-Abelian Lie group whose Lie algebra is semisimple, i.e., it is a direct sum of simple Lie algebras} Lie groups \cite{hall2013lie} like $\mathrm{SU}(n)\ (n\geq2), \mathrm{SO}(n)\ (n \geq 3), \mathrm{SL}(n) \text{ and } \mathrm{Sp}(n)\ (n \geq 1)$, it holds the total antisymmetry of the structure constants $c_{ki}^j p_j p^k = 0$. Consequently, for this wide range of matrix Lie groups, the Lie derivative operator \eqref{eq:kinetic_vf} of the kinetic part simplifies so that explicit forms of the Lie derivative operators read
\begin{equation}\label{eq:Lie_derivative_explicit}
    \hat{\V} = -\e_i(\V) \frac{\partial}{\partial p_i} \quad \text{and} \quad \hat{\T} = p^i \e_i,
\end{equation}
so that Hamilton's equations become
\begin{align}\label{eq:EoM_LieGroup}
\begin{alignedat}{2}
    \dot{\bfP} &= - \e_i(\V) \frac{\partial \bfP}{\partial p_i} = - \e_i(\V) \bfT^i, & \bfP(0) &= \bfP_0, \\
    \dot{\bfQ} &= p^i \e_i(\bfQ) = - p^i \bfT_i \bfQ = - \bfP \bfQ,\quad & \bfQ(0) &= \bfQ_0.
\end{alignedat}
\end{align}
The formal solution of the subsystems 
\begin{align*}
    \begin{pmatrix}
    \dot{\bfP} \\ \dot{\bfQ}
\end{pmatrix} &= \begin{pmatrix}
    \hat{\V} \bfP \\ \0
\end{pmatrix} = \begin{pmatrix}
    -\e_i(\V) \T^i \\ \0
\end{pmatrix}, & \quad \begin{pmatrix}
    \dot{\bfP} \\ \dot{\bfQ}
\end{pmatrix} &= \begin{pmatrix}
    \0 \\ \hat{\T} \bfQ
\end{pmatrix} = \begin{pmatrix}
    \0 \\ -\bfP \bfQ
\end{pmatrix} 
\intertext{can be expressed in terms of the matrix exponential $\exp(t \A) = \sum\nolimits_{k=0}^\infty (t\A)^k/k!$ and the Lie derivative operators via}
\varphi_t^{\{1\}}(\bfP_0,\bfQ_0) &= \begin{pmatrix}
    \exp(t \hat{\V})\bfP_0 \\ \bfQ_0
\end{pmatrix} = \begin{pmatrix}
    \bfP_0 - t \e_i(\V(\bfQ_0))\bfT^i \\
    \bfQ_0
\end{pmatrix}, & \varphi_t^{\{2\}}(\bfP_0,\bfQ_0) &= \begin{pmatrix}
    \bfP_0 \\ \exp(t \hat{\T})\bfQ_0
\end{pmatrix} = \begin{pmatrix}
    \bfP_0 \\ \exp(-t \bfP_0) \bfQ_0
\end{pmatrix},
\end{align*}
for the potential and the kinetic part, respectively.

{\remark[System updates]{When having many components $(\p_\ell,\q_\ell) \in T^*\G,\ \ell=1,\ldots,L$, one can trivially generalize the fundamental two-form, resulting in sums over all components $\bomega = - \sum\nolimits_{\ell=1}^L \diff \p_\ell$. Analogously, the locally acting kinetic energy generalizes to $\T(\p) = \sum\nolimits_{\ell=1}^L \tfrac{1}{2} \langle \p_\ell,\p_\ell \rangle$. Consequently, Hamilton's equations become 
\begin{align*}
\begin{alignedat}{2}
    \dot{\bfP}_\ell &= - \e_i(\V)_\ell \bfT^i, \quad & \bfP_\ell(0) &= \bfP_{\ell,0}, \\
    \dot{\bfQ}_\ell &= - \bfP_\ell \bfQ_\ell, \quad & \bfQ_\ell(0) &= \bfQ_{\ell,0}, \; \ell = 1,\ldots,L.
\end{alignedat}
\end{align*}}}
\remark[Classical mechanics]{In classical Hamiltonian mechanics, the phase space is $(\p,\q) \in \mathbb{R}^d \times \mathbb{R}^d$. With the fundamental two-form $\bomega = \diff\q \wedge \diff\p$, one obtains the well-known Hamiltonian equations of motion $$\begin{pmatrix}
    \dot{\p} \\ \dot{\q}
\end{pmatrix} = \begin{pmatrix}
    - \frac{\partial \Ham}{\partial \q} \\ \frac{\partial \Ham}{\partial \p}
\end{pmatrix} = \begin{pmatrix}
    - \V_{\q}(\q) \\ \T_{\p}(\p)
\end{pmatrix}, \quad \begin{pmatrix}
    \p(0) \\ \q(0)
\end{pmatrix} = \begin{pmatrix}
    \p_0 \\ \q_0
\end{pmatrix}.$$
}\normalfont
For the sake of clarity and coherence, we will keep the notation \eqref{eq:EoM_LieGroup} for Hamiltonian mechanics on matrix Lie groups. 
However, it is important to emphasize that the results in this paper also hold for system updates, as well as for classical mechanics in Euclidean space.

\subsection{Geometric integration}
The Hamiltonian flow $\varphi_t(\p_0,\q_0)$ is characterized by the following properties. 
\begin{itemize}
    \item \textbf{Energy conservation:} The Hamiltonian is an invariant of the flow, i.e., \begin{equation*}
        \frac{\diff}{\diff t} \Ham(\varphi_t(\p_0,\q_0)) = 0.
    \end{equation*}
    \item \textbf{Time-reversibility:} The Hamiltonian flow is time-reversible, i.e., \begin{equation*}
        \rho \circ \varphi_t \circ \rho \circ \varphi_t(\p_0,\q_0) = (\p_0,\q_0), 
        \quad \rho(\p,\q) = (-\p,\q).
    \end{equation*}
    \item \textbf{Symplecticity:} The Hamiltonian flow is symplectic, i.e., $\diff\bomega = 0$.
    A direct consequence of the symplecticity of the Hamiltonian flow is the preservation of volume,
    \begin{equation*}
        \left\lvert \mathrm{det} \frac{\partial \varphi_t(\p_0,\q_0)}{\partial (\p_0,\q_0)} \right\rvert = 1.
    \end{equation*} 
    \item \textbf{Closure property:} The Hamiltonian flow satisfies $\varphi_t(\p_0,\q_0) \in T^*\G$ for all $t > 0$ provided that $(\p_0,\q_0) \in T^*\G$.
\end{itemize}
We demand the numerical scheme $\Phi_h(\p_0,\q_0) = (\p_1, \q_1) \approx \varphi_h(\p_0,\q_0)$ to preserve the time-reversibility, volume-preservation, as well as the closure property of the Hamiltonian flow. For \emph{time-reversibility}, one gets the criterion \begin{equation}\label{eq:time-reversibility}
    \rho \circ \Phi_h \circ \rho \circ \Phi_h(\p_0,\q_0) = (\p_0,\q_0).
\end{equation}
The numerical integration scheme is \emph{volume-preserving} if it holds 
\begin{equation}\label{eq:volume-preservation}
    \left\lvert \det \frac{\partial \Phi_h(\p_0,\q_0)}{\partial(\p_0,\q_0)} \right\rvert = 1.
\end{equation}
The \emph{closure property} demands
\begin{equation}\label{eq:closure_property}
    \Phi_h(\p_0,\q_0) = (\p_1,\q_1) \in T^*\G, 
\end{equation}
provided $(\p_0,\q_0) \in T^*\G$.

\section{Hessian-free force-gradient integrators}\label{sec:Hessian-free_FGI}
An efficient approach for geometric numerical integration of Hamiltonian systems of the form \eqref{eq:Hamiltonian} is given by force-gradient integrators (FGIs) \cite{omelyan2003symplectic,kennedy2009force}. 
In this section, we will introduce the new class of Hessian-free FGIs. Before introducing the adaption, we will start with a brief recapitulation of FGIs based on their introduction in \cite{omelyan2003symplectic}.

\subsection{Force-gradient integrators}
As we have seen in the previous section, it is possible to compute the flows of the subsystems using the exponential map and the Lie derivative operators \eqref{eq:Lie_derivative_explicit}. Hence it is possible to compute an approximation to Hamilton's equations \eqref{eq:EoM_LieGroup} by composing evaluations of these exact flows, resulting in splitting methods \cite{mclachlan2002splitting}
\begin{align*}%\label{eq:splitting_method}
    \Psi^{(0)} = \mathrm{Id}, \quad \Psi^{(j)} = \Psi^{(j-1)} \circ \varphi_{a_j h}^{\{2\}} \circ \varphi_{b_j h}^{\{1\}}, \quad j=1,\ldots,P.
\end{align*}
Since we demand the numerical integration scheme to preserve the time-reversibility \eqref{eq:time-reversibility}, it has to hold either $a_1 = 0$ (velocity version) or $b_P = 0$ (position version). The remaining composition of $2P-1$ exponentials has to be symmetric, i.e., it has to hold $a_{j+1} = a_{P-j+1}$ and $b_j = b_{P-j+1}$ for the velocity version, $a_j = a_{P-j+1}$ and $b_j = b_{P-j}$ for the position version. Applying the Baker--Campbell--Hausdorff (BCH) formula, the overall integrator $\Psi_h := \Psi_h^{(P)}$ can be written as an exponential \cite{omelyan2003symplectic}
\begin{subequations}\label{eq:BCH}
\begin{equation}\label{eq:BCH_exponential}
    \exp\left((\nu \hat{\T} + \sigma \hat{\V})h + \mathcal{O}_3 h^3 + \mathcal{O}_5 h^5 + \mathcal{O}_7 h^7 + \mathcal{O}(h^9)\right)
\end{equation}
with $\nu \coloneqq \sum_{j=1}^P a_j$, $\sigma \coloneqq \sum_{j=1}^P b_j$, and
\begin{align}
    \mathcal{O}_3 &= \alpha \left[\hat{\T},\left[\hat{\T},\hat{\V}\right]\right] + \beta \left[\hat{\V},\left[\hat{\T}, \hat{\V}\right]\right], \label{eq:BCH_order3}\\
\begin{split}\label{eq:BCH_order5}
    \mathcal{O}_5 &= \gamma_1 \left[\hat{\T},\left[\hat{\T},\left[\hat{\T},\left[\hat{\T}, \hat{\V}\right]\right]\right]\right]  +\gamma_2 \left[\hat{\T},\left[\hat{\T},\left[\hat{\V},\left[\hat{\T}, \hat{\V}\right]\right]\right]\right] \\
    &\quad + \gamma_3 \left[\hat{\V},\left[\hat{\T},\left[\hat{\T},\left[\hat{\T}, \hat{\V}\right]\right]\right]\right] + \gamma_4 \left[\hat{\V},\left[\hat{\V},\left[\hat{\T},\left[\hat{\T}, \hat{\V}\right]\right]\right]\right],\footnotemark 
\end{split}\\
\begin{split}\label{eq:BCH_order7}
    \mathcal{O}_7 &= \zeta_1 \left[\hat{\V},\left[\hat{\V},\left[\hat{\T},\left[\hat{\V},\left[\hat{\T},\left[\hat{\V},\hat{\T}\right]\right]\right]\right]\right]\right] + \zeta_2 \left[\hat{\V},\left[\hat{\V},\left[\hat{\V},\left[\hat{\T},\left[\hat{\T},\left[\hat{\V},\hat{\T}\right]\right]\right]\right]\right]\right] \\
    &\quad + \zeta_3 \left[\hat{\V},\left[\hat{\V},\left[\hat{\T},\left[\hat{\T},\left[\hat{\T},\left[\hat{\V},\hat{\T}\right]\right]\right]\right]\right]\right] + \zeta_4 \left[\hat{\V},\left[\hat{\T},\left[\hat{\V},\left[\hat{\T},\left[\hat{\T},\left[\hat{\V},\hat{\T}\right]\right]\right]\right]\right]\right] \\
    &\quad + \zeta_5 \left[\hat{\T},\left[\hat{\V},\left[\hat{\V},\left[\hat{\T},\left[\hat{\T},\left[\hat{\V},\hat{\T}\right]\right]\right]\right]\right]\right] + \zeta_6 \left[\hat{\T},\left[\hat{\V},\left[\hat{\T},\left[\hat{\V},\left[\hat{\T},\left[\hat{\V},\hat{\T}\right]\right]\right]\right]\right]\right]\\
    &\quad + \zeta_7 \left[\hat{\V},\left[\hat{\T},\left[\hat{\T},\left[\hat{\T},\left[\hat{\T},\left[\hat{\V},\hat{\T}\right]\right]\right]\right]\right]\right] + \zeta_8 \left[\hat{\T},\left[\hat{\V},\left[\hat{\T},\left[\hat{\T},\left[\hat{\T},\left[\hat{\V},\hat{\T}\right]\right]\right]\right]\right]\right]\\
    &\quad + \zeta_9 \left[\hat{\T},\left[\hat{\T},\left[\hat{\V},\left[\hat{\T},\left[\hat{\T},\left[\hat{\V},\hat{\T}\right]\right]\right]\right]\right]\right] + \zeta_{10} \left[\hat{\T},\left[\hat{\T},\left[\hat{\T},\left[\hat{\T},\left[\hat{\T},\left[\hat{\V},\hat{\T}\right]\right]\right]\right]\right]\right],
\end{split}
\end{align}
\end{subequations}
where $[\cdot,\cdot]$ denotes the commutator $[\A,\B] \coloneqq \A\B - \B\A$.\footnotetext{neglecting two zero-valued commutators  $[\hat{\V},[\hat{\V},[\hat{\V},[\hat{\T},\hat{\V}]]]]$ and $[\hat{\T},[\hat{\V},[\hat{\V},[\hat{\T},\hat{\V}]]]]$} 
To obtain splitting methods of convergence order $p>2$, one has to cancel the order-three error term $\mathcal{O}_3$, consisting of the two commutators $[\hat{\T},[\hat{\T},\hat{\V}]]$ and $[\hat{\V},[\hat{\T},\hat{\V}]]$. Without using negative time steps, that are used for example in composition techniques \cite{yoshida1990construction,suzuki1990fractal,omelyan2002construction}, it is not possible to cancel both commutators at once \cite{yoshida1990construction}. For Hamiltonian systems of the form \eqref{eq:Hamiltonian}, the second commutator $\hat{\C} \coloneqq [\hat{\V},[\hat{\T},\hat{\V}]]$ exhibits a special structure
\begin{align}\label{eq:force-gradient_term}
    \hat{\C} &= 2 \hat{\V} \hat{\T} \hat{\V} 
    = 2 \e^j(\V) \e_j\e_i(\V)\frac{\partial}{\partial p_i}
\end{align}
that only depends on the generalized coordinates $q_j$, and is called \emph{force-gradient term} (FG-term). Extending the momentum updates $\varphi_{b_j h}^{\{1\}}$ by including evaluations of the force-gradient term, $\exp(b_j \hat{\V} + c_j \hat{\C})$, results in a FGI. This allows to choose the coefficients $a_j, b_j$ such that the first commutator $[\hat{\T},[\hat{\T},\hat{\V}]]$ vanishes and then the coefficients $c_j$ can be chosen to remove the remaining order-3 error term $[\hat{\V},[\hat{\T},\hat{\V}]]$, resulting in higher-order integrators. A symmetric FGI consists of $2P-1$ exponentials that are either momentum updates 
\begin{align}\label{eq:FGI_mom-update}
    \exp\left(b^{(n)} h \hat{\V} + c^{(n)} h^3 \hat{\C}\right)
\end{align}
or position updates 
\begin{align}\label{eq:FGI_pos-update}
    \exp\left(a^{(n)} h \hat{\T}\right)
\end{align}
with coefficients $a^{(n)},b^{(n)},c^{(n)}$ that are related to the coefficients $a_j,b_j,c_j$ via 
\begin{align*}
    j = \begin{cases}
        \tfrac{P+2}{2} - n, & P \text{\ even (velocity version),}\\
        \tfrac{P+1}{2}+n, & P \text{\ odd (velocity version),}\\
        \tfrac{P}{2}+n, & P \text{\ even (position version),}\\
        \tfrac{P+1}{2}-n, & P \text{\ odd (position version),}
    \end{cases}
\end{align*}
for $n=1,\ldots,\lfloor \tfrac{P}{2}\rfloor$. 
A FGI is constructed as follows. Starting from a central single-exponential operator 
\begin{subequations}\label{eq:force-gradient_integrator}
\begin{equation}\label{eq:FGI_central-exp}
    \Psi^{(0)} = \begin{cases} 
    \exp( a_{(P+2)/2}h \hat{T} ) , & P \text{ even (velocity version)},  \\
    \exp( b_{(P+1)/2}h \hat{\V} + c_{(P+1)/2}h^3 \hat{\C} ) , &  P \text{ odd (velocity version)}, \\
    \exp( b_{P/2}h \hat{\V} + c_{P/2}h^3 \hat{\C} ) , &  P \text{ even (position version)}, \\
    \exp( a_{(P+1)/2}h \hat{\T} ) , &  P \text{ odd (position version)},
\end{cases}
\end{equation}
the integrator is obtained by consecutively applying the following symmetric transformations
\begin{align}
\begin{split}
    \Psi^{(n)} &= \exp\left(b^{(n)} h \hat{\V} + c^{(n)} h^3 \hat{\C}\right) \circ \exp\left(a^{(n)} h \hat{\T}\right) \circ \Psi^{(n-1)} \circ \exp\left(a^{(n)} h \hat{\T}\right) \circ \exp\left(b^{(n)} h \hat{\V} + c^{(n)} h^3 \hat{\C}\right), \label{eq:FGI_velocity-version}
\end{split}
\intertext{for $P$ odd (velocity version) and $P$ even (position version), and}
\begin{split}
    \Psi^{(n)} &= \exp\left(a^{(n)} h \hat{\T}\right) \circ \exp\left(b^{(n)} h \hat{\V} + c^{(n)} h^3 \hat{\C}\right) \circ \Psi^{(n-1)} \circ \exp\left(b^{(n)} h \hat{\V} + c^{(n)} h^3 \hat{\C}\right) \circ \exp\left(a^{(n)} h \hat{\T}\right), \label{eq:FGI_position-version}
\end{split}
\end{align}
\end{subequations}
for $P$ even (velocity version) and $P$ odd (position version), $n=1,\ldots,\lfloor \tfrac{P}{2}\rfloor$. By applying the BCH formula, the overall FGI $\Psi_h := \Psi_h^{(\lfloor P/2 \rfloor)}$ again takes the form \eqref{eq:BCH}. Thanks to the recursive definition \eqref{eq:force-gradient_integrator} of the integrator, the multipliers $\nu,\sigma,\alpha,\beta$, $\gamma_1,\ldots,\gamma_4$, and $\zeta_1,\ldots,\zeta_{10}$ can be determined using recursive formulations stated in \cite{omelyan2003symplectic}.
\medbreak
\noindent\textbf{Order conditions.} The FGI \eqref{eq:force-gradient_integrator} has convergence order $p$ if its representation as an exponential \eqref{eq:BCH} satisfies 
\begin{align*}
    \exp\!\left((\nu \hat{\T} + \sigma \hat{\V})h \!+\! \mathcal{O}_3 h^3 + \mathcal{O}_5 h^5 + \mathcal{O}_7 h^7 + \mathcal{O}(h^9)\right) = \exp\!\left(h(\hat{\T} \!+\! \hat{\V})\right) \!+\! \mathcal{O}(h^{p+1}).
\end{align*}
Hence the order conditions up to order seven are given by the recursive formulations of the multipliers $\nu,\sigma,\alpha,\beta,$  $\gamma_1,\ldots,\gamma_4$, and $\zeta_1,\ldots,\zeta_{10}$ in \cite{omelyan2003symplectic}.
\medbreak
\noindent\textbf{Geometric integration.} The momentum updates \eqref{eq:FGI_mom-update} of the FGI with $c^{(n)} \neq 0$ can be regarded as an evaluation of the exact flow, corresponding to a modified potential 
\begin{equation*}%\label{eq:modified_potential_FGI}
    \V_{\mathrm{FG}}(\q) = \V(\q) - \tfrac{c^{(n)} h^2}{b^{(n)}} \V_{\hat{\C}}(\q),
\end{equation*} 
where $\V_{\hat{\C}}(\q)$ is solved by $\e_i(\V_{\hat{\C}}) = 2 \e^j(\V) \e_j\e_i(\V)$, that is again Hamiltonian. Due to the symmetric construction of \eqref{eq:force-gradient_integrator}, the integrator is time-reversible \eqref{eq:time-reversibility}. 
As the integrator is a composition of exact flows that are symplectic, the overall integration scheme is symplectic and thus volume-preserving \eqref{eq:volume-preservation}. The integrator also satisfies the closure property \eqref{eq:closure_property} as a) the position updates \eqref{eq:FGI_pos-update} remain unchanged and b) the momentum updates \eqref{eq:FGI_mom-update} still yield momenta $\bfP \in \frakg$.
\medbreak
\noindent\textbf{Backward error analysis.} It is a well-known quantity of symplectic integration schemes that they preserve a nearby \emph{shadow Hamiltonian} $\tilde{\Ham}$ exactly \cite{kennedy2013shadow}. Given a symplectic integrator of order $p$, it holds $\tilde{\Ham} = \Ham + \mathcal{O}(h^{p})$. As the FGI can be written as an exponential \eqref{eq:BCH}, it is straight-forward to determine an explicit expression for $\tilde{\Ham}$. Replacing the commutators of the vector fields $\hat{\V}, \hat{\T}$ by Poisson brackets of the zero-forms $\V,\T$ gives the shadow Hamiltonian 
\begin{align}\label{eq:Shadow_Hamiltonian_FGI}
\begin{split}
    \tilde{\Ham}_{\mathrm{FG}} = \nu \T &+ \sigma \V + h^2 \left(\alpha \left\{ \T, \left\{ \T,\V \right\} \right\} + \beta \left\{ \V, \left\{ \T,\V \right\} \right\}  \right) \\
    &+ h^4 \left( \gamma_1 \left\{\T,\left\{\T,\left\{\T,\left\{\T, \V\right\}\right\}\right\}\right\} +\gamma_2 \left\{\T,\left\{\T,\left\{\V,\left\{\T, \V\right\}\right\}\right\}\right\} \right. \\
    &\qquad \left. + \gamma_3 \left\{\V,\left\{\T,\left\{\T,\left\{\T, \V\right\}\right\}\right\}\right\} + \gamma_4 \left\{\V,\left\{\V,\left\{\T,\left\{\T, \V\right\}\right\}\right\}\right\} \right) \\ 
    &+ h^6 \left( \zeta_1 \left\{\V,\left\{\V,\left\{\T,\left\{\V,\left\{\T,\left\{\V,\T\right\}\right\}\right\}\right\}\right\}\right\} + \zeta_2 \left\{\V,\left\{\V,\left\{\V,\left\{\T,\left\{\T,\left\{\V,\T\right\}\right\}\right\}\right\}\right\}\right\} \right. \\
    &\quad + \zeta_3 \left\{\V,\left\{\V,\left\{\T,\left\{\T,\left\{\T,\left\{\V,\T\right\}\right\}\right\}\right\}\right\}\right\} + \zeta_4 \left\{\V,\left\{\T,\left\{\V,\left\{\T,\left\{\T,\left\{\V,\T\right\}\right\}\right\}\right\}\right\}\right\} \\
    &\quad + \zeta_5 \left\{\T,\left\{\V,\left\{\V,\left\{\T,\left\{\T,\left\{\V,\T\right\}\right\}\right\}\right\}\right\}\right\} + \zeta_6 \left\{\T,\left\{\V,\left\{\T,\left\{\V,\left\{\T,\left\{\V,\T\right\}\right\}\right\}\right\}\right\}\right\}\\
    &\quad + \zeta_7 \left\{\V,\left\{\T,\left\{\T,\left\{\T,\left\{\T,\left\{\V,\T\right\}\right\}\right\}\right\}\right\}\right\} + \zeta_8 \left\{\T,\left\{\V,\left\{\T,\left\{\T,\left\{\T,\left\{\V,\T\right\}\right\}\right\}\right\}\right\}\right\}\\
    &\quad \left. + \zeta_9 \left\{\T,\left\{\T,\left\{\V,\left\{\T,\left\{\T,\left\{\V,\T\right\}\right\}\right\}\right\}\right\}\right\} + \zeta_{10} \left\{\T,\left\{\T,\left\{\T,\left\{\T,\left\{\T,\left\{\V,\T\right\}\right\}\right\}\right\}\right\}\right\}\right) + \mathcal{O}(h^8)
\end{split}
\end{align}
that is exactly preserved by the FGI \eqref{eq:force-gradient_integrator}.

\example[fourth-order FGI]{\label{ex:FGI4} One example of a fourth-order FGI that demands only five stages and a single evaluation of the FG-term is given by
\begin{equation}\label{eq:FGI4}
    \exp\left(\tfrac{h}{6} \hat{\V}\right) \exp\left(\tfrac{h}{2} \hat{\T}\right) \exp\left(\tfrac{2h}{3} \hat{\V} + \tfrac{h^3}{72} \hat{\C}\right) \exp\left(\tfrac{h}{2} \hat{\T}\right) \exp\left(\tfrac{h}{6} \hat{\V}\right).
\end{equation}
This integrator exactly preserves the shadow Hamiltonian
\begin{align*}%\label{eq:shadow_FGI4}
    \tilde{\Ham}_{\mathrm{FGI4}} &= \T + \V + h^4 \left( \tfrac{1}{2880} \left\{\T,\left\{\T,\left\{\T,\left\{\T, \V\right\}\right\}\right\}\right\} + \tfrac{1}{2880} \left\{\T,\left\{\T,\left\{\V,\left\{\T, \V\right\}\right\}\right\}\right\} \right. \\
    &\qquad \left. + \tfrac{1}{2160} \left\{\V,\left\{\T,\left\{\T,\left\{\T, \V\right\}\right\}\right\}\right\} + \tfrac{1}{4320} \left\{\V,\left\{\V,\left\{\T,\left\{\T, \V\right\}\right\}\right\}\right\} \right) + \mathcal{O}(h^6).     
\end{align*}
}\normalfont

In general, one has to derive the FG-term $\hat{\C}$ for each particular Hamiltonian system by computing the Hessian of the potential. An approximation of the FG-term proposed in \cite{wisdom1996symplectic} enables the use of the force-gradient integrator \eqref{eq:FGI4} by replacing the FG-term by an additional force evaluation. We will generalize this idea to the entire class of FGIs.

\subsection{Approximation of the force-gradient step}
The derivation of the FG-term \eqref{eq:force-gradient_term} can be quite tedious. In \cite{shcherbakov2017adapted}, for example, the FG-term for the two-dimensional Schwinger model has been computed. In other applications, e.g.\ the application for Hybrid Monte Carlo (HMC) \cite{duane1987hybrid} simulations in lattice QCD on four-dimensional lattices \cite{knechtli2017lattice}, it becomes even more complicated to derive the FG-term. It is possible to approximate the force-gradient step \eqref{eq:FGI_mom-update} by a composed evaluation of the force $\e_i(\V)$ as it has been shown for Example \ref{ex:FGI4} in \cite{yin2011improving}.   
 A generalization to this approximation for arbitrary FGIs reads
\begin{align}\label{eq:Yin_derivation}
\begin{split}
    \exp\!\left(b^{(n)} h \hat{\V} + c^{(n)} h^3 \hat{\C}\right) &=
    \exp\!\left(\!-b^{(n)} h \e_i(\V) \frac{\partial}{\partial p_i} + 2 c^{(n)} h^3 \e^j(\V)  \e_j\e_i(\V) \frac{\partial}{\partial p_i} \!\right), \\
    &= \exp\!\left(-b^{(n)} h \left(\mathrm{Id} - \frac{2c^{(n)} h^2}{b^{(n)}} \e^j(\V) \e_j \right) \e_i(\V) \frac{\partial}{\partial p_i}\right), \\ 
    &= \exp\!\left(-b^{(n)} h \exp\left(- \frac{2c^{(n)} h^2}{b^{(n)}} F^j \e_j \right) \e_i(\V) \frac{\partial}{\partial p_i} \right) + \mathcal{O}(h^5).
\end{split}
\end{align}
Here, $F^j \e_j = \e^j(\V)_{\q} \e_j$ denotes that the vector field is frozen at $(\p,\q)$ \cite{CrouchGrossman} so that $\e_j$ acting on $F^j$ is defined to be zero.
Moreover, performing a Taylor expansion verifies the equality 
$$ \exp\left(-\frac{2c^{(n)}h^2}{b^{(n)}} F^j \e_j \right) \e_i(\V)(\bfQ) = \e_i(\V)\left( \exp\left( - \frac{2 c^{(n)}h^2}{b^{(n)}} F^j \bfT_j \right) \bfQ \right).$$

Consequently, this approximation allows one to replace the momentum update \eqref{eq:FGI_mom-update} by the following two-step procedure:
\begin{enumerate}
    \item Compute a temporary position update via $\bfQ' = \exp\left(- \frac{2c^{(n)} h^2 }{b^{(n)}} F^j \bfT_j\right) \bfQ.$
    \item Compute a usual momentum update of step size $b^{(n)} h$ by using the temporary position update for evaluating the force term $-\e_i(\V)$, i.e., $\bfP_{b^{(n)} h} = \bfP_0 - b^{(n)} h \e_i(\V)(\bfQ') \bfT^i.$ 
\end{enumerate}

Although this two-step procedure is used in implementations, the following equivalent expression is of importance for performing the error analysis.

\theorem[Hessian-free force-gradient step]{The two-step procedure is equivalent to replacing the momentum update \eqref{eq:FGI_mom-update} by  \begin{equation}\label{eq:FGI_exponential_Yin}
    \exp\left(b^{(n)} h \hat{\D}(h,b^{(n)},c^{(n)})\right) := \exp\left(b^{(n)} h \sum\limits_{k=0}^\infty 
    \frac{2^k c^{(n)^k} h^{2k}}{k! \cdot b^{(n)^k}}
    \Big(\mathcal{\widetilde{V}}\hat{\T}\Big)^k \hat{\V}\right),
\end{equation}
with $\mathcal{\widetilde{V}} := - \e_i(\V)_{\q} \frac{\partial}{\partial p_i}$.
\begin{proof}
When applying the exponential \eqref{eq:FGI_mom-update} to $(\p,\q)$, one obtains 
\begin{align*}\exp\left( b^{(n)} h \hat{\V} + c^{(n)} h^3 \hat{\C} \right) \begin{pmatrix}
    \p \\ \q
\end{pmatrix} &= \exp\left( b^{(n)} h \hat{\V} + 2 c^{(n)} h^3 \hat{\V} \hat{\T} \hat{\V} \right) \begin{pmatrix}
    \p \\ \q
\end{pmatrix}
= \exp\left( b^{(n)} h \left( \mathrm{Id} + \tfrac{2c^{(n)}h^2}{b^{(n)}} \hat{\V}\hat{\T} \right) \hat{\V}\right) \begin{pmatrix}
    \p \\ \q
\end{pmatrix} \\
&= \exp\left( b^{(n)} h \left( \mathrm{Id} + \tfrac{2c^{(n)}h^2}{b^{(n)}} \mathcal{\widetilde{V}} \hat{\T} \right) \hat{\V} \right)
\begin{pmatrix}
    \p \\ \q
\end{pmatrix} = \exp\left( b^{(n)} h \exp\left( \tfrac{2c^{(n)}h^2}{b^{(n)}} \mathcal{\widetilde{V}} \hat{\T}\right) \hat{\V} \right) \begin{pmatrix}
    \p \\ \q
\end{pmatrix} + \mathcal{O}(h^5).
\end{align*}
Applying the definition of the matrix exponential to the inner exponential concludes the proof. \endproof
\end{proof}}\normalfont

Consequently, the symmetric transformations \eqref{eq:FGI_velocity-version} and \eqref{eq:FGI_position-version} become
\begin{subequations}\label{eq:FGI_mom-update_Yin}
\begin{align}
\begin{split}\label{eq:Hf-FGI_velocity}
    \Psi^{(n)} &= \exp\left(b^{(n)} h \hat{\D}(h,b^{(n)},c^{(n)})\right) \circ \exp\left(a^{(n)} h \hat{\T}\right) \circ \Psi^{(n-1)} \circ \exp\left(a^{(n)} h \hat{\T}\right) \circ \exp\left(b^{(n)} h \hat{\D}(h,b^{(n)},c^{(n)})\right),
\end{split}\\
\begin{split}\label{eq:Hf-FGI_position}
    \Psi^{(n)} &= \exp\left(a^{(n)} h \hat{\T}\right) \circ \exp\left(b^{(n)} h \hat{\D}(h,b^{(n)},c^{(n)})\right) \circ \Psi^{(n-1)} \circ \exp\left(b^{(n)} h \hat{\D}(h,b^{(n)},c^{(n)})\right) \circ \exp\left(a^{(n)} h \hat{\T}\right),
\end{split}
\end{align}
\end{subequations}
respectively. We refer to FGIs that utilize the approximation \eqref{eq:Yin_derivation} as \emph{Hessian-free FGIs}. From now on, we will refer to FGIs evaluating the FG-term~\eqref{eq:force-gradient_term} as \textit{exact FGIs}.
\medbreak
\noindent\textbf{Order conditions.} As the approximation \eqref{eq:FGI_exponential_Yin} introduces an error of $\mathcal{O}(h^5)$, the order conditions up to order four do not change.
Based on the expression \eqref{eq:FGI_exponential_Yin}, it is possible to explicitly state the formulae of the error terms for Hessian-free FGIs. Of particular interest are the changes in the error terms \eqref{eq:BCH_order5} and \eqref{eq:BCH_order7}. The new error terms read
\begin{align}\label{eq:new_error_terms}
\begin{split}
    \mathcal{O}'_5 &= \mathcal{O}_5 + \gamma_5\  \mathcal{\widetilde{V}}\hat{\T}\mathcal{\widetilde{V}}\hat{\T}\hat{\V}, \\
    \mathcal{O}'_7 &= \mathcal{O}_7 + \zeta_{11} \mathcal{\widetilde{V}}\hat{\T}\mathcal{\widetilde{V}}\hat{\T}\mathcal{\widetilde{V}} \hat{\T} \hat{\V} + \zeta_{12} \hat{\T}\mathcal{\widetilde{V}}\hat{\T}\mathcal{\widetilde{V}}\hat{\T}\hat{\V}\hat{\T} + \zeta_{13} \hat{\T}\hat{\T}\mathcal{\widetilde{V}}\hat{\T}\mathcal{\widetilde{V}}\hat{\T}\hat{\V},
\end{split}
\end{align}
where the new multipliers $\gamma_5,\zeta_{11},\zeta_{12},\zeta_{13}$ can be computed based on recursive relations. Starting from
\begin{subequations}\label{eq:gamma5}
\begin{align*}
    \gamma_5^{(0)} &:= \begin{cases} 
    0, & P \text{ even (velocity version)},  \\
    2 c_{(P+1)/2}^2/b_{(P+1)/2}, & P \text{ odd (velocity version)}, \\
    2 c_{P/2}^2/b_{P/2}, & P \text{ even (position version)}, \\
    0, &  P \text{ odd (position version)},
\end{cases} & 
\zeta_{11}^{(0)} &:= \begin{cases} 
    0, & P \text{ even (velocity version)},  \\
    4 c_{(P+1)/2}^3/(3 b_{(P+1)/2}^2), &  P \text{ odd (velocity version)}, \\
    4 c_{P/2}^3/(3 b_{P/2}^2), &  P \text{ even (position version)}, \\
    0, &  P \text{ odd (position version)},
\end{cases}
\end{align*}
and $\zeta_{12}^{(0)} = \zeta_{13}^{(0)} = 0$, the recursive relations for the symmetric transformations \eqref{eq:Hf-FGI_velocity} become
\begin{align*}
    \gamma_5^{(n)} &= \gamma_5^{(n-1)} + 4 c^{(n)^2}/b^{(n)}, &
    \zeta_{11}^{(n)} &= \zeta_{11}^{(n-1)} + \left(8 c^{(n)^3}/b^{(n)} + 2 \sigma^{(n)} \nu^{(n)} c^{(n)^2}\right)/(3 b^{(n)}), \\
    \zeta_{12}^{(n)} &= \zeta_{12}^{(n-1)} + a^{(n)^2} \gamma_5^{(n-1)}/3 - 2 \nu^{(n)^2} c^{(n)^2} / (3 b^{(n)}), &
    \zeta_{13}^{(n)} &= \zeta_{13}^{(n-1)} - a^{(n)^2} \gamma_5^{(n-1)}/6 + \nu^{(n)^2} c^{(n)^2} / (3 b^{(n)}),
\intertext{and for the symmetric transformations \eqref{eq:Hf-FGI_position}}
    \gamma_5^{(n)} &= \gamma_5^{(n-1)} + 4 c^{(n)^2}/b^{(n)}, &
    \zeta_{11}^{(n)} &= \zeta_{11}^{(n-1)} + \left(8 c^{(n)^3}/b^{(n)} + 2 \sigma^{(n)} \nu^{(n)} c^{(n)^2}\right)/(3 b^{(n)}), \\
    \zeta_{12}^{(n)} &= \zeta_{12}^{(n-1)} + a^{(n)^2} \gamma_5^{(n)}/3 - 2 \nu^{(n)^2} c^{(n)^2} / (3 b^{(n)}), &
    \zeta_{13}^{(n)} &= \zeta_{13}^{(n-1)} - a^{(n)^2} \gamma_5^{(n)}/6 + \nu^{(n)^2} c^{(n)^2} / (3 b^{(n)}),
\end{align*}
\end{subequations}
for $n=1,\ldots,\lfloor \tfrac{P}{2} \rfloor$.
\medbreak
\noindent\textbf{Geometric integration.} The Hessian-free force-gradient step \eqref{eq:FGI_exponential_Yin} yields updated momenta $\bfP \in \frakg$, whereas the position updates \eqref{eq:FGI_pos-update} remain unchanged. Consequently, the closure property \eqref{eq:closure_property} is still satisfied. The approximation does not affect the time-reversibility of the integrator. Furthermore, Hessian-free FGIs are a composition of shears that are volume-preserving so that the overall integrator preserves the volume. Hessian-free FGIs are no longer symplectic. The modified momentum updates \eqref{eq:FGI_exponential_Yin} are only symplectic if the Hessian of the potential commutes with itself when evaluated at (in general, different) $\bfQ$. Generally, this is not the case.
\medbreak
\noindent\textbf{Backward error analysis.} Due to the additional error terms, the modified differential equation is not Hamiltonian, i.e.\ it does not preserve a nearby shadow Hamiltonian. Instead, the modified differential equation of the Hessian-free FGI reads 
\begin{align}\label{eq:modified_system}
    \begin{pmatrix}
        \dot{\bfP} \\ \dot{\bfQ} 
    \end{pmatrix} &= \begin{pmatrix}
        \{ \tilde{\Ham}_{\mathrm{FG}}, \bfP \} \\ \{\tilde{\Ham}_{\mathrm{FG}}, \bfQ\} 
    \end{pmatrix} + \begin{pmatrix}
        h^4 \gamma_5 \mathcal{\widetilde{V}} \hat{\T} \mathcal{\widetilde{V}} \hat{\T} \hat{\V} \bfP + h^6 \zeta_{11} \mathcal{\widetilde{V}}\hat{\T}\mathcal{\widetilde{V}}\hat{\T}\mathcal{\widetilde{V}} \hat{\T} \hat{\V} \bfP + h^6 \zeta_{13} \hat{\T}\hat{\T}\mathcal{\widetilde{V}}\hat{\T}\mathcal{\widetilde{V}}\hat{\T}\hat{\V} \bfP \\ h^6 \zeta_{12} \hat{\T}\mathcal{\widetilde{V}}\hat{\T}\mathcal{\widetilde{V}}\hat{\T}\hat{\V}\hat{\T}\bfQ 
    \end{pmatrix} + \mathcal{O}(h^8),
\end{align}
where the first term corresponds to the shadow Hamiltonian \eqref{eq:Shadow_Hamiltonian_FGI} of the exact FGI using the same set of coefficients.
The remaining terms introduce an energy drift. In general, without making any further assumption on the potential $\V$ of the Hessian-free FGI, the energy will have a linear drift of size $\mathcal{O}(th^{\max\{4,p\}})$ for a Hessian-free FGI of order $p$ and a time interval of length $t$.
Note that one may choose the coefficients of the Hessian-free FGI to cancel drift terms so that one ends up with Hessian-free FGIs of order $p$ where the energy will have a linear drift of size $\mathcal{O}(t h^r)$ with $r > p$. 

\example[Hessian-free fourth-order FGI]{\label{ex:adapted_FGI4}Based on the fourth-order exact FGI \eqref{eq:FGI4}, one gets the Hessian-free FGI 
\begin{align}\label{eq:FGI4_Yin}
    \exp\left(\tfrac{h}{6} \hat{\V}\right) \exp\left(\tfrac{h}{2} \hat{\T}\right) \exp\left(\tfrac{2h}{3} \hat{\D}(h,\tfrac{2}{3},\tfrac{1}{72})\right) \exp\left(\tfrac{h}{2} \hat{\T}\right) \exp\left(\tfrac{h}{6} \hat{\V}\right)
\end{align}
of convergence order four. With $\gamma_5 = 1/1728 \neq 0$, the numerical energy has a linear drift of size $\mathcal{O}(th^4)$.
}\normalfont 
Despite the lack of symplecticity, Hessian-free FGIs are interesting integrators and are particularly useful in the following situations.
\begin{enumerate}
    \item[a)] The evaluation of the FG-term $\hat{\C}$ is significantly more expensive than an evaluation of the force $\e_i(\V)$,
    \item[b)] An analytical expression for the FG-term $\hat{\C}$ is not available, e.g.\ because its derivation is tedious.
\end{enumerate}
In certain applications, such as particles moving in gravitational or coulomb fields, evaluations of the FG-term are less expensive than force evaluations due to the specific structure of the potential. Here, Hessian-free FGIs are less efficient and not of practical relevance unless the exact FG-term is not available. In case of expensive evaluations (e.g. in molecular dynamics simulations), it demands a refined analysis of Hessian-free FGIs to clarify whether Hessian-free FGIs allow for a more efficient computational process. For this purpose, the next section provides a classification and derivation of Hessian-free FGIs with up to 11 stages, similar to the investigations made in \cite{omelyan2003symplectic}. 

\section{Derivation of Hessian-free force-gradient integrators}\label{sec:integrator_derivation}
In this section, we will provide a complete classification of Hessian-free FGIs with up to eleven stages. The coefficients are chosen so that the maximum possible convergence order $p$ is achieved. If the solution is not unique, the free parameters can be optimized by minimizing a (weighted) norm of the multipliers in $\mathcal{O}'_{p+1}$. By expanding the commutators in \eqref{eq:new_error_terms}, one ends up with the following contributions containing the new error terms:
\begin{itemize}
    \item $4 \gamma_4 \hat{\V} \hat{\T} \hat{\V} \hat{\T} \hat{\V} + \gamma_5 \mathcal{\widetilde{V}} \hat{\T} \mathcal{\widetilde{V}} \hat{\T} \hat{\V}$,
    \item $-8 \zeta_1 \hat{\V}\hat{\T}\hat{\V}\hat{\T}\hat{\V}\hat{\T}\hat{\V} + \zeta_{11} \mathcal{\widetilde{V}}\hat{\T}\mathcal{\widetilde{V}}\hat{\T}\mathcal{\widetilde{V}}\hat{\T}\hat{\V}$,
    \item $8 \zeta_6 \hat{\T}\hat{\V}\hat{\T}\hat{\V}\hat{\T}\hat{\V}\hat{\T} + \zeta_{12} \hat{\T}\mathcal{\widetilde{V}}\hat{\T}\mathcal{\widetilde{V}}\hat{\T}\hat{\V}\hat{\T} $,
    \item $(-3 \zeta_4 - 4 \zeta_6)\hat{\T}\hat{\T}\hat{\V}\hat{\T}\hat{\V}\hat{\T}\hat{\V} + \zeta_{13} \hat{\T}\hat{\T}\mathcal{\widetilde{V}}\hat{\T}\mathcal{\widetilde{V}}\hat{\T}\hat{\V}$.
\end{itemize}
Whereas we choose to minimize the norm 
\begin{align}\label{eq:Err3}
    \mathrm{Err}_3 \coloneqq \sqrt{\alpha^2 + \beta^2}
\end{align}
in case of $p=2$, we use weighted norms for $p=4$ and $p=6$ incorporating the aforementioned relations. For $p=4$, we choose the norm 
\begin{align}\label{eq:Err5}
    \mathrm{Err}_5 &\coloneqq \sqrt{\gamma_1^2 + \gamma_2^2 + \gamma_3^2 + \gamma_4^2 + \left(\tfrac{1}{4}\gamma_5\right)^2},
\end{align} 
which, in case of non-gradient schemes ($\gamma_5 = 0$), results in the usual Euclidean norm used in \cite{omelyan2003symplectic}. The weighted norm incorporates that the contribution of $\gamma_5$ is four times less than $\gamma_4$. Analogously, we choose for $p=6$ 
\begin{align}\label{eq:Err7}
    \mathrm{Err}_7 \coloneqq \sqrt{\sum\limits_{j=1}^{10} \zeta_j^2 + \left(\tfrac{1}{8} \zeta_{11}\right)^2 + \left(\tfrac{1}{8} \zeta_{12}\right)^2 + \left(\tfrac{7}{24} \zeta_{13}\right)^2}.
\end{align}
Here, the weight of $7/24$ for $\zeta_{13}$ is derived by taking the mean value of $1/3$ and $1/4$, so that for non-gradient schemes ($\zeta_{11} = \zeta_{12} = \zeta_{13} = 0$) the weighted norm again reduces to the norm used in \cite{omelyan2003symplectic}.
Moreover, the efficiency of the integrators has been measured using the formula
\begin{equation}\label{eq:efficiency_measure}
    \mathrm{Eff}^{(p)} = \frac{1}{(n_f)^p \cdot \mathrm{Err}_{p+1}}.
\end{equation}
Note that in order to reduce the number of unknowns, we will take into account in advance the symmetry of the FGI, as well as the order-1 conditions $\sum\nolimits_{j=1}^P a_j = \sum\nolimits_{j=1}^P b_j = 1$ throughout this section. For any number of stages $s$, we obtain different variants based on the coefficients $c_j$ set to zero in order to reduce the number $n_f$ of force evaluations. All variants are summarized by the total number of force evaluations $n_f$ per time step, the achieved convergence order $p$, and the number of degrees of freedom $v$, followed by the set of coefficients that is minimizing the (weighted) norm of the leading error term $\mathrm{Err}_{p+1}$. In addition, the corresponding efficiency according to \eqref{eq:efficiency_measure} is stated. Note that for velocity algorithms, the number of force recalculations $n_f$ is reduced by one since the last modified momentum update of time step $n$ and the first one of time step $n+1$ coincide so that the force does not have to be re-evaluated.

\subsection{Hessian-free force-gradient integrators with s=3 stages (P=2)}
For $P=2$, it is only possible to achieve convergence order $p=2$. The free parameters are chosen so that the norm \eqref{eq:Err3}
is minimized.
\subsubsection{Velocity versions} 
The velocity version of a convergent three-stage Hessian-free FGI reads 
\begin{align*}%\label{eq:aFGI_s3v}
    \Psi_h = \exp\left(\tfrac{h}{2} \hat{\D}(h,1/2,c_1) \right) \exp\left( h \hat{\T} \right) \exp\left(\tfrac{h}{2} \hat{\D}(h,1/2,c_1) \right)
\end{align*}
with free parameter $c_1$. The method has convergence order two with leading error terms $\alpha = 1/12$ and $\beta = 1/24 + 2 c_1$. 

\noindent\textit{Variant 1 ($n_f = 2$, $p=2$, $v=1$).} 
\begin{align}\label{eq:DAD}
\begin{split}
    c_1 &= -1/48, \qquad
    \mathrm{Err}_3 = 1/12, \qquad
    \mathrm{Eff} = 3.
\end{split}
\end{align} 

\noindent\textit{Variant 2 ($n_f = 1$, $p=2$, $v=0$).} For $c_1 = 0$, one obtains the well-known velocity version of the Störmer--Verlet scheme \cite{HairerLubichWanner}, resulting in
\begin{align}\label{eq:BAB} 
\begin{split}
    \mathrm{Err}_3 = \sqrt{5}/24 \approx 0.0932, \quad \text{and} \quad 
    \mathrm{Eff} = 24/\sqrt{5} \approx 10.73.
\end{split}
\end{align}

\subsubsection{Position versions} 
The position version for $P=2$ reads 
\begin{equation*}%\label{eq:aFGI_s3p}
    \Psi_h = \exp\left(\tfrac{h}{2} \hat{\T}\right)\exp\left(h \hat{\D}(h,1,c_1)\right)\exp\left(\tfrac{h}{2} \hat{\T}\right)
\end{equation*}
with free parameter $c_1$. The scheme has convergence order two and leading error terms $\alpha = -1/24$ and $\beta = -1/12 + c_1$.

\noindent\textit{Variant 1 ($n_f = 2$, $p=2$, $v=1$).} 
\begin{align}\label{eq:ADA}
    c_1 = 1/12, \qquad \mathrm{Err}_3 = 1/24, \qquad \mathrm{Eff} = 6.
\end{align} 

\noindent\textit{Variant 2 ($n_f = 1$, $p=2$, $v=0$).} Setting $c_1 = 0$ yields the position version of the Störmer--Verlet method \cite{HairerLubichWanner} with
\begin{align}\label{eq:ABA}
    \mathrm{Err}_3 = \sqrt{5}/24\approx 0.0932, \quad \text{and} \quad 
    \mathrm{Eff} = 24/\sqrt{5} \approx 10.7.
\end{align}

\subsubsection{Remarks} 
As all schemes for $s=3$ cannot achieve a higher convergence order than two, the optimal values for $c_1$, as well as the principal error term do not change compared to the results in \cite{omelyan2003symplectic}. However, the efficiency of the Hessian-free FGIs with $c_1 \neq 0$ increases since in \cite{omelyan2003symplectic} evaluating the FG-term is assumed to be twice as expensive as the force.

\subsection{Hessian-free force-gradient integrators with s=5 stages (P=3)}
Considering $P=3$, it is possible to obtain Hessian-free FGIs of convergence order $p=4$. For fourth-order methods, we minimize $\mathrm{Err}_5$ defined in \eqref{eq:Err5}
and denote its minimum by $\gamma_{\min}$.

\subsubsection{Velocity versions}
The velocity version
\begin{align*}%\label{eq:aFGI_s5v}
    \Psi_h &= \exp\left(b_1 h \hat{\D}(h,b_1,c_1)\right) \exp\left(\tfrac{h}{2} \hat{\T}\right) \exp\left((1-2b_1) h \hat{\D}(h,1-2b_1,c_2)\right) \exp\left(\tfrac{h}{2} \hat{\T}\right) \exp\left(b_1 h \hat{\D}(h,b_1,c_1)\right)
\end{align*}
has three free parameters $b_1,c_1,c_2$. Solving the order-three conditions while keeping $c_1$ as a free parameter gives 
\begin{align}\label{eq:aFGI_s5v_eqs}
    b_1 = 1/6, \qquad c_2 = 1/72 - 2c_1.  
\end{align}

\noindent\textit{Variant 1 ($n_f = 4$, $p=4$, $v=1$).} The minimization of $\mathrm{Err}_5$ yields
\begin{align}\label{eq:DADAD}
\begin{split}
    c_1 = -0.000881991367333, \quad 
    c_2 =   0.015652871623554, \quad
    \gamma_{\min} \approx 0.000625, \quad \mathrm{Eff} \approx 6.25.
\end{split}
\end{align}

\noindent\textit{Variant 2 ($n_f = 3$, $p=4$, $v=0$).} Setting $c_1 = 0$, \eqref{eq:aFGI_s5v_eqs} yields the unique solution
\begin{align}\label{eq:BADAB}
\begin{split}
    c_2 = 1/72, \quad \gamma_{\min} \approx 0.000728, \quad \mathrm{Eff} \approx 16.96.
\end{split}
\end{align}

\noindent\textit{Variant 3 ($n_f = 3$, $p=4$, $v=0$).} For $c_2 = 0$, \eqref{eq:aFGI_s5v_eqs} yields  
\begin{align}\label{eq:DABAD}
\begin{split}
    c_1 = 1/144, \quad \gamma_{\min} \approx 0.00335, \quad \mathrm{Eff} \approx 3.68.
\end{split}
\end{align}

\noindent\textit{Variant 4 ($n_f = 2$, $p=2$, $v=1$).} For $c_1 = c_2 = 0$, we obtain a non-gradient algorithm, i.e., the approximation of the FG-term does not apply. Hence we obtain exactly the same coefficients and principal error term as in \cite{omelyan2003symplectic}, namely
\begin{align}\label{eq:BABAB}
    b_1 &= \frac{1}{2} - \frac{(2 \sqrt{326} + 36)^{1/3}}{12} + \frac{1}{6(2 \sqrt{326} + 36)^{1/3}} \approx 0.1931833275037836, &
    \mathrm{Err}_3 &\approx 0.00855, & \mathrm{Eff} &\approx 29.24. 
\end{align} 

\subsubsection{Position versions}
The position version 
\begin{align*}%\label{eq:aFGI_s5p}
    \Psi_h &= \exp\left(a_1 h \hat{\T} \right) \exp\left(\tfrac{h}{2} \hat{\D}(h,1/2,c_1) \right) \exp\left((1-2a_1)h \hat{\T} \right) \exp\left(\tfrac{h}{2} \hat{\D}(h,1/2,c_1) \right) \exp\left(a_1 h \hat{\T} \right) 
\end{align*}
has two degrees of freedom, namely $a_1,c_1$. 

\noindent\textit{Variant 1 ($n_f = 4$, $p=4$, $v=0$).} We obtain two real solutions where the optimal one reads
\begin{align}\label{eq:ADADA}
    a_1 = \frac{1}{2}\left(1 - \frac{1}{\sqrt{3}} \right),\quad c_1 = \frac{1}{48}(2 - \sqrt{3}), \quad
    \gamma_{\min} \approx 0.000718, \quad \mathrm{Eff} \approx 5.44.
\end{align}

\noindent\textit{Variant 2 ($n_f = 2$, $p=2$, $v=1$).} Setting $c_1 = 0$, the Hessian-free FGI reduces to a non-gradient scheme that does not make use of the approximation to the FG-term $\hat{\C}$. The global minimum of $\mathrm{Err}_3$ is achieved at
\begin{align}\label{eq:ABABA}
    a_1 &= \frac{1}{2} - \frac{(2 \sqrt{326} + 36)^{1/3}}{12} + \frac{1}{6(2 \sqrt{326} + 36)^{1/3}} \approx 0.1931833275037836, &
    \mathrm{Err}_3 &\approx 0.00855, & \mathrm{Eff} &\approx 29.24.
\end{align}

\subsubsection{Remarks} All variants except \eqref{eq:DADAD} have the same optimal solutions as the exact FGIs in \cite{omelyan2003symplectic}. However, their principal error term and efficiency change due to the additional error term $\gamma_5$ (except for the non-gradient versions). In case of \eqref{eq:DADAD}, the additional error term also results in a difference in the coefficients, particularly in $c_1$ and thus in $c_2$. The Hessian-free FGI \eqref{eq:BADAB} coincides with \eqref{eq:FGI4_Yin} which was the first Hessian-free FGI derived in \cite{yin2011improving} and has been applied to the Schwinger model \cite{shcherbakov2017adapted} in a nested integrator. It is the most efficient Hessian-free FGI with $s=5$ stages. The non-gradient schemes \eqref{eq:BABAB} and \eqref{eq:ABABA} define the most efficient integrators of order two within this framework.

\subsection{Hessian-free force-gradient integrators with s=7 stages (P=4)}
For $P=4$, one can again derive Hessian-free FGIs up to convergence order four. However, the additional number of stages allows to reduce the principal error term $\mathrm{Err}_5$. Furthermore, it is now possible to derive non-gradient schemes of order $p=4$.
\subsubsection{Velocity versions}
The velocity version with $P=4$ reads 
\begin{align*}
    \Psi_h &= \exp\left( b_1 h \hat{\D}(h,b_1,c_1) \right) \exp\left( a_2 h \hat{\T} \right) \exp\left( (\tfrac{1}{2} - b_1)h \hat{\D}(h,\tfrac{1}{2} - b_1,c_2)  \right) \exp\left( (1-2a_2)h \hat{\T}\right) \\
    & \qquad \times  \exp\left( (\tfrac{1}{2} - b_1)h \hat{\D}(h,\tfrac{1}{2} - b_1,c_2)  \right) \exp\left( a_2 h \hat{\T} \right) \exp\left( hb_1 \hat{\D}(h,b_1,c_1) \right)
\end{align*}
where $a_2,b_1,c_1,c_2$ are free parameters. The order-three conditions have the solutions 
\begin{align}\label{eq:aFGI_s7v_eqs}
    b_1 = \frac{1}{12} \left( 6 + \frac{1}{a_2(a_2 - 1)} \right), \quad
    c_1 = -\frac{1}{288} \left( 6 + 288c_2 - \frac{1}{a_2(a_2 - 1)^2}\right).
\end{align}

\noindent\textit{Variant 1 ($n_f = 6$, $p=4$, $v=2$).} Minimizing $\mathrm{Err}_5$ with respect to the two free parameters $a_2,c_2$ results in 
\begin{equation}
\begin{alignedat}{3}\label{eq:DADADAD}
    a_2 &= 0.273005515864808, \qquad &
    b_1 &= 0.080128674198082, \qquad &
    c_1 &= 0.000271601364672, \\ 
    c_2 &= 0.002959399979707, &
    \gamma_{\min} &\approx 0.0000275, & \mathrm{Eff} &\approx 28.09. 
\end{alignedat}
\end{equation}

\noindent\textit{Variant 2 ($n_f = 5$, $p=4$, $v=1$).} At $c_1 = 0$, the global minimum is reached at
\begin{equation}\label{eq:BADADAB}
  \begin{alignedat}{3}
    a_2 &= 0.281473422092232, \qquad & b_1&= 0.087960811032557,\qquad &
    c_2 &= 0.003060423791562, \\ 
    \gamma_{\min}&\approx 0.0000498, & \mathrm{Eff} &\approx 32.12. & &
  \end{alignedat}
\end{equation}

\noindent\textit{Variant 3 ($n_f = 4$, $p=4$, $v=1$).} Setting $c_2 = 0$, one finds
\begin{equation}\label{eq:DABABAD}
\begin{alignedat}{3}
    a_2 &=  0.258529167713908, \qquad &
    b_1 &=  0.065274481323251, \qquad &
    c_1 &=  0.003595899064589,\\ 
    \gamma_{\min} &\approx 0.000891, & \mathrm{Eff} &\approx 4.38. & & 
\end{alignedat}
\end{equation}

\noindent\textit{Variant 4 ($n_f = 3$, $p=4$, $v=0$).} For $c_1 = 0$ and $c_2 = 0$, one obtains the velocity version of the Forest-Ruth algorithm \cite{forest1992sixth} with
\begin{align}\label{eq:BABABAB}
    a_1 &= 1/(2 - 2^{1/3}), &
    b_1 &= 1/(2\cdot(2-2^{1/3})), &
    \gamma_{\min} &\approx 0.0383, & \mathrm{Eff} &\approx 0.32.
\end{align}

\subsubsection{Position versions} The position counterpart has the form
\begin{align*}%\label{eq:aFGI_s7p}
%\begin{split}
    \Psi_h &= \exp\left(a_1 h \hat{\T} \right) \exp\left( b_1 h \hat{\D}(h,b_1,c_1) \right) \exp\left((\tfrac{1}{2} - a_1)h \hat{\T} \right) \exp\left( (1-2b_1)h \hat{\D}(h,1-2b_1,c_2) \right) \\
    & \qquad \times \exp\left((\tfrac{1}{2} - a_1)h \hat{\T} \right) \exp\left( b_1 h \hat{\D}(h,b_1,c_1) \right) \exp\left(a_1 h \hat{\T} \right),
%\end{split}    
\end{align*}
where $a_1,b_1,c_1,c_2$ are free parameters. Analogously to the velocity counterpart, one can solve the third-order conditions for $a_1$ and $c_1$, leading to 
\begin{align}\label{eq:aFGI_s7p_eqs}
\begin{split}
    a_1 = \frac{1}{2} \pm \frac{1}{\sqrt{24 b_1}}, \quad
    c_1 = \frac{1}{24} \left( 1 - 12 c_2 \pm \sqrt{6b_1} (1-b_1) \right).
\end{split}    
\end{align}

\noindent\textit{Variant 1 ($n_f = 6$, $p=4$, $v=2$).} One reaches the minimum of $\mathrm{Err}_5$ at 
\begin{equation}\label{eq:ADADADA}
\begin{alignedat}{3}
    a_1 &= 0.116438749543126,\qquad &
    b_1 &= 0.283216992495952,\qquad &
    c_1 &= 0.001247201195115, \\
    c_2 &= 0.002974030329635, &
    \gamma_{\min} &\approx 0.0000200, & 
    \mathrm{Eff} &\approx 38.57.
\end{alignedat}
\end{equation}

\noindent\textit{Variant 2 ($n_f = 5$, $p=4$, $v=1$).} Putting $c_2 = 0$, one achieves the minimum at
\begin{equation}\label{eq:ADABADA}
\begin{alignedat}{3}
    a_1 &= 0.136458051118946, \qquad &
    b_1 &=  0.315267858070664, \qquad &
    c_1 &= 0.002427032834125, \\
    \gamma_{\min} &\approx 0.0000844, & 
    \mathrm{Eff} &\approx 18.95. & & 
\end{alignedat}
\end{equation}

\noindent\textit{Variant 3 ($n_f = 4$, $p=4$, $v=1$).} Setting $c_1 = 0$, the minimum is achieved at
\begin{equation}\label{eq:ABADABA}
\begin{alignedat}{3}
    a_1 &= 0.089775972994422, \qquad &
    b_1 &= 0.247597680043986, \qquad &
    c_2 &= 0.006911440413815, \\
    \gamma_{\min} &\approx 0.000149, &
    \mathrm{Eff} &\approx 26.19. & &
\end{alignedat}
\end{equation}

\noindent\textit{Variant 4 ($n_f = 3$, $p=4$, $v=0$).} The non-gradient version, obtained by setting $c_1 = c_2 = 0$, is the position version of the algorithm by Forest and Ruth \cite{forest1992sixth} which is given by
\begin{align}\label{eq:ABABABA}
    a_1 &= 1/(2\cdot(2-2^{1/3})), &
    b_1 &= 1/(2 - 2^{1/3}), &
    \gamma_{\min} &\approx 0.0283, & \mathrm{Eff} &\approx 0.44.
\end{align}

\subsubsection{Remarks} The most efficient seven-stage algorithm is given by the position version \eqref{eq:ADADADA} with an efficiency of approximately $38.57$. Compared to the most efficient five-stage Hessian-free FGI of order four, namely the velocity version \eqref{eq:BADAB} with $\mathrm{Eff} \approx 16.96$. Consequently, using the same number of total force recalculations, the seven-stage algorithm \eqref{eq:ADADADA} will reduce the global errors approximately by a factor of $38.57/16.96 \approx 2.27$ compared to the five-stage Hessian-free FGI \eqref{eq:BADAB}. The non-gradient versions \eqref{eq:BABABAB} and \eqref{eq:ABABABA} now define methods of order four. These methods are equivalent to applying Yoshida's triple-jump composition \cite{yoshida1990construction} to \eqref{eq:BAB} and \eqref{eq:ABA}, respectively. The principal error term of the non-gradient schemes is relatively large, resulting in less efficient algorithms. In practice, one should either choose force-gradient schemes (e.g., the five-stage Hessian-free FGI \eqref{eq:BADAB} reduces the global errors approximately by a factor of $16.96/0.44 \approx 38.55$ while using the same number of total force evaluations) or non-gradient schemes with more stages which will be presented in the proceeding sections.

\subsection{Hessian-free force-gradient integrators with s=9 stages (P=5)}
Taking $P=5$ enables the derivation of a Hessian-free FGI of convergence order $p=6$. In this case, we aim for minimizing the norm \eqref{eq:Err7}
and denote its minimum by $\zeta_{\min}$,
otherwise we remain at minimizing $\mathrm{Err}_5$ according to \eqref{eq:Err5}. 

\subsubsection{Velocity versions}
The velocity version of a nine-stage Hessian-free FGI takes the form 
\begin{align*}
    \Psi_h &= \exp\left(b_1 h \hat{\D}(h,b_1,c_1) \right) \exp\left(a_2 h \hat{\T}\right) \exp\left(b_2 h \hat{\D}(h,b_2,c_2) \right) \\
    &\qquad \times \exp\left( (\tfrac{1}{2} - a_2)h \hat{\T}\right) \exp\left((1-2(b_1 + b_2))h \hat{\D}(h,1-2(b_1+b_2),c_3) \right) \\
    &\qquad \times \exp\left( (\tfrac{1}{2} - a_2)h \hat{\T}\right) \exp\left(b_2 h \hat{\D}(h,b_2,c_2) \right)\exp\left(a_2 h \hat{\T}\right)\exp\left(b_1 h \hat{\D}(h,b_1,c_1) \right)
\end{align*} 
and is characterized by six degrees of freedom $a_2,b_1,b_2,c_1,c_2,c_3$.

\noindent\textit{Variant 1 ($n_f = 7$, $p=6$, $v=0$).} Despite there is a new order condition due to $\gamma_5$, we find the same unique solution as in \cite{omelyan2003symplectic}, namely
\begin{align}\label{eq:BADADADAB}
\begin{split}
    a_2 &= \frac{1}{2} + \frac{\sqrt[3]{675 + 75 \sqrt{6}}}{30} + \frac{5}{2 \sqrt[3]{675 + 75\sqrt{6}}}, \quad 
    b_1 = \frac{a_2}{3},\quad
    b_2 = - \frac{5 a_2}{3}(a_2 - 1), \quad
    c_1 = 0, \\
    c_2 &= - \frac{5 a_2^2}{144} + \frac{a_2}{36} - \frac{1}{288}, \quad
    c_3 = \frac{1}{144} - \frac{a_2}{36} \left( \frac{a_2}{2} + 1 \right), \quad 
    \zeta_{\min} \approx 0.00154, \quad 
    \mathrm{Eff} \approx 0.0055.
\end{split}
\end{align}
Since the unique solution satisfies $c_1 = 0$, this variant coincides with the variant of demanding $c_1 = 0$.

\noindent\textit{Variant 2 ($n_f = 7$, $p=4$, $v=3$).} For $c_3 = 0$, minimizing the norm of the leading error term results in
\begin{equation}\label{eq:DADABADAD}
\begin{alignedat}{3}
    a_2 &= 0.227758000273404, \qquad &
    b_1 &= 0.070935378258660, \qquad &
    b_2 &= 0.322911610232109, \\ 
    c_1 &=  0.000067752132787, &
    c_2 &= 0.001597508440746, & 
    \gamma_{\min} &\approx 0.0000101,\; \mathrm{Eff} \approx 41.06.
\end{alignedat}
\end{equation}

\noindent\textit{Variant 3 ($n_f = 6$, $p=4$, $v=3$).} By demanding $c_2 = 0$, one gets 
\begin{equation}\label{eq:DABADABAD}
\begin{alignedat}{3}
    a_2 &= 0.197279141794602, \qquad &
    b_1  &= 0.060885008530668, \qquad &
    b_2 &= 0.288579639891554, \\
    c_1  &= 0.000429756946246, &
    c_3 &= 0.002373498029145, & 
    \gamma_{\min} &\approx 0.0000130,\; \mathrm{Eff} \approx 59.33.
\end{alignedat}
\end{equation}

\noindent\textit{Variant 4 ($n_f = 6$, $p=4$, $v=2$).} At $c_1 = 0$ and $c_3 = 0$, the minimization of $\mathrm{Err}_5$ yields 
\begin{equation}\label{eq:BADABADAB}
\begin{alignedat}{3}
    a_2 &= 0.219039425103133,\qquad &
    b_1  &= 0.068466565514186, \qquad &
    b_2 &= 0.311000565033563, \\
    c_2  &= 0.001602470431500, &
    \gamma_{\min} &\approx 0.0000105, & \mathrm{Eff} &\approx 73.45.
\end{alignedat}
\end{equation}

\noindent\textit{Variant 5 ($n_f = 5$, $p=4$, $v=2$).} For $c_1 = 0$ and $c_2 = 0$, we achieve 
\begin{equation}\label{eq:BABADABAB}
\begin{alignedat}{3}
    a_2 &= 0.200395293638238,\qquad & 
    b_1  &= 0.073943321445602, \qquad &
    b_2 &= 0.258244950046509, \\
    c_3  &= 0.003147048491590, &
    \gamma_{\min} &\approx 0.0000651, & \mathrm{Eff} &\approx 24.57.
\end{alignedat}
\end{equation}

\noindent\textit{Variant 6 ($n_f = 5$, $p=4$, $v=2$).} Putting $c_2 = 0$ and $c_3 = 0$, one finds the minimum
\begin{equation}\label{eq:DABABABAD}
\begin{alignedat}{3}
    a_2 &= 0.190585159174513,\qquad &
    b_1  &= 0.036356798097337, \qquad &
    b_2 &= 0.340278911234329, \\
    c_1  &= 0.002005691094612, &
    \gamma_{\min} &\approx 0.000336, & \mathrm{Eff} &\approx 4.76. 
\end{alignedat}
\end{equation}

\noindent\textit{Variant 7 ($n_f = 4$, $p=4$, $v=1$).} Finally, choosing $c_1 = 0$, $c_2 = 0$, $c_3 = 0$ and minimizing $\mathrm{Err}_5$ results in the non-gradient scheme \cite{omelyan2002construction}
\begin{equation}\label{eq:BABABABAB}
\begin{alignedat}{3}
    a_2 &= 0.520943339103990, \qquad &
    b_1 &= 0.164498651557576,  \qquad &
    b_2 &= 1.235692651138917, \\
    \gamma_{\min} &\approx 0.000654, & 
    \mathrm{Eff} &\approx 5.97. & &
\end{alignedat}
\end{equation}

\subsubsection{Position versions}
The position version 
\begin{align*}
    \Psi_h &= \exp\left(a_1 h \hat{\T} \right) \exp\left(b_1 h \hat{\D}(h,b_1,c_1) \right) \exp\left(a_2 h \hat{\T}\right)\exp\left((\tfrac{1}{2}-b_1)h \hat{\D}(h,\tfrac{1}{2}-b_1,c_2) \right) \exp\left((1-2(a_1+a_2))h \hat{\T} \right) \\
    &\qquad \times \exp\left((\tfrac{1}{2}-b_1)h \hat{\D}(h,\tfrac{1}{2}-b_1,c_2) \right) \exp\left(a_2 h \hat{\T}\right) \exp\left(b_1 h \hat{\D}(h,b_1,c_1) \right) \exp\left(a_1 h \hat{\T} \right)
\end{align*}
comes with five degrees of freedom $a_1,a_2,b_1,c_1,c_2$, only allowing for the derivation of Hessian-free FGIs up to convergence order $p=4$. The order conditions up to order four yield the equations 
\begin{align}\label{eq:aFGI_s9p_eqs}
\begin{split}
    a_1 &= \left(3 - 6a_2 + 12a_2 b_1 \mp \sqrt{3} \sqrt{1 - 24 a_2^2 b_1 + 48 a_2^2 b_1^2}\right)/6, \\
    c_2 &= \left(2-12 a_2 b_1 + 24 a_2 b_1^2 \mp \sqrt{3} \sqrt{1 - 24a_2^2 b_1 + 48 a_2^2 b_1^2} - 48 c_1\right)/48,
\end{split}
\end{align}
so that three degrees of freedom $a_2,b_1,c_1$ are remaining. 

\noindent\textit{Variant 1 ($n_f = 8$, $p=4$, $v=3$).} The norm of the leading error term is minimized at sign minus with 
\begin{equation}\label{eq:ADADADADA}
\begin{alignedat}{3}
    a_1 &= 0.094471605659163, \qquad &
    a_2 &= 0.281057227947299, \qquad &
    b_1 &= 0.227712700174579, \\
    c_1 &=  0.000577062053569, &
    c_2 &= 0.000817399268485, & \gamma_{\min} &\approx 0.00000501,\; \mathrm{Eff} \approx 48.71.
\end{alignedat}
\end{equation}

\noindent\textit{Variant 2 ($n_f = 6$, $p=4$, $v=2$).} Setting $c_1 = 0$ results in a minimum achieved at sign minus, given by 
\begin{equation}\label{eq:ABADADABA}
\begin{alignedat}{3}
    a_1 &= 0.047802682977081, \qquad &
    a_2 &=  0.265994592108478, \qquad &
    b_1 &= 0.143282503449494, \\
    c_2 &=  0.002065558490728, &
    \gamma_{\min}&\approx 0.0000346, & \mathrm{Eff} &\approx 22.32.
\end{alignedat}
\end{equation}

\noindent\textit{Variant 3 ($n_f = 6$, $p=4$, $v=2$).} For $c_2 = 0$, we obtain (at sign minus)
\begin{equation}\label{eq:ADABABADA}
\begin{alignedat}{3}
    a_1 &= 0.118030603246046, \qquad &
    a_2 &=  0.295446189611111, \qquad &
    b_1 &= 0.273985556386628, \\
    c_1 &=  0.001466561305710, &
    \gamma_{\min}&\approx 0.0000471, & \mathrm{Eff} &\approx 16.39.
\end{alignedat}
\end{equation}

\noindent\textit{Variant 4 ($n_f = 4$, $p=4$, $v=1$).} With $c_1 = 0$ and $c_2 = 0$, one finds the non-gradient scheme 
\begin{equation}\label{eq:ABABABABA}
\begin{alignedat}{3}
    a_1 &=  0.178617895844809,\qquad &
    a_2 &= -0.066264582669818, \qquad &
    b_1 &=  0.712341831062606, \\
    \gamma_{\min} &\approx 0.000610, &
    \mathrm{Eff} &\approx 6.40. & &
\end{alignedat}
\end{equation}

\subsubsection{Remarks} For $s=9$, we obtain the first (and at the same time unique) Hessian-free FGI of convergence order $p=6$, namely the velocity version \eqref{eq:BADADADAB}. It is the same unique solution as in \cite{omelyan2003symplectic}, although one has to satisfy the new order condition due to $\gamma_5$. Since some variants of eleven-stage Hessian-free FGIs will be affected by order reduction, this can be regarded as a lucky coincidence. All other variants define order-four integrators. The most efficient one is given by the velocity version \eqref{eq:BADABADAB} with an efficiency of $\mathrm{Eff} \approx 73.45$. Compared to the most efficient seven-stage variant \eqref{eq:ADADADA}, it allows to reduce the global errors approximately by a factor of $73.45/38.57 \approx 1.9$ while keeping the total number of force recalculations fixed.
  
\subsection{Hessian-free force-gradient integrators with s=11 stages (P=6)}
We finally consider the case of Hessian-free FGIs with $P=6$. For the cases $P\leq 5$, each variant of Hessian-free FGIs has the same convergence order as the common FGIs in \cite{omelyan2003symplectic}, i.e., the additional order conditions did not result in an order reduction. For $P=6$, however, some variants are affected by order reduction.

\subsubsection{Velocity versions}
In velocity form, the eleven-stage propagation takes the form
\begin{align*}
    \Psi_h &= \exp\left( b_1 h \hat{\D}(h,b_1,c_1) \right) \exp\left(a_2 h \hat{\T} \right) \exp\left( b_2 h \hat{\D}(h,b_2,c_2) \right) \exp\left(a_3 h \hat{\T} \right) \exp\left( (\tfrac{1}{2} \!-\! (b_1\!+\!b_2))h \hat{\D}(h,\tfrac{1}{2} \!-\! (b_1\!+\!b_2),c_3) \right) \\
    &\qquad \times \exp\left( (1-2(a_2+a_3))h \hat{\T} \right) \exp\left( (\tfrac{1}{2} \!-\! (b_1\!+\!b_2))h \hat{\D}(h,\tfrac{1}{2} \!-\! (b_1\!+\!b_2),c_3) \right) \exp\left(a_3 h \hat{\T} \right) \exp\left( b_2 h \hat{\D}(h,b_2,c_2) \right) \\
    &\qquad \times \exp\left(a_2 h \hat{\T} \right) \exp\left( b_1 h \hat{\D}(h,b_1,c_1) \right)
\end{align*}
with seven degrees of freedom $a_2,a_3,b_1,b_2,c_1,c_2,c_3$.

\noindent\textit{Variant 1 ($n_f = 9$, $p=6$, $v=0$).} The only real solution is given by 
\begin{equation}\label{eq:BADADADADAB}
\begin{alignedat}{3}
    a_2 &= 0.270990466773838,\qquad & 
    a_3 &= 0.635374358266882, \qquad &
    b_1 &= 0.090330155591279, \\
    b_2 &= 0.430978044876253, &
    c_1 &= 0, & 
    c_2 &= 0.002637435980472, \\
    c_3 &= -0.000586445610932, & 
    \zeta_{\min} &\approx 0.00000699, & \mathrm{Eff} &\approx 0.27.
\end{alignedat}
\end{equation}
Note that \eqref{eq:BADADADADAB} also covers the variant demanding $c_1=0$.

\noindent\textit{Variant 2 ($n_f = 8$, $p=4$, $v=4$).} For $c_2 = 0$, the best solution is achieved at
\begin{equation}\label{eq:DABADADABAD}
\begin{alignedat}{3}
    a_2 &= 0.068597474282941,\qquad & 
    a_3 &= 0.284851197274498, \qquad &
    b_1 &= -0.029456704762871, \\
    b_2 &= 0.228751459942521, &
    c_1 &= 0.000410146066173, & 
    c_3 &= 0.001249935251564, \\
    \gamma_{\min} &\approx 0.00000355, & \mathrm{Eff} &\approx 68.84. & &
\end{alignedat}
\end{equation}

\noindent\textit{Variant 3 ($n_f = 8$, $p=4$, $v=4$).} Putting $c_3 = 0$ yields 
\begin{equation}\label{eq:DADABABADAD}
\begin{alignedat}{3}
    a_2 &= 0.203263079324187,\qquad & 
    a_3 &= 0.200698071607808, \qquad &
    b_1 &= 0.066202529912271, \\
    b_2 &= 0.267856111220228, &
    c_1 &= 0.000012570620797, & 
    c_2 &= 0.001042408779514, \\
    \gamma_{\min} &\approx 0.00000519, & \mathrm{Eff} &\approx 47.08. & &
\end{alignedat}
\end{equation}

\noindent\textit{Variant 4 ($n_f = 7$, $p=4$, $v=3$).} Setting $c_1 = 0$ and $c_2 = 0$, the optimal solution looks as
\begin{equation}\label{eq:BABADADABAB}
\begin{alignedat}{3}
    a_2 &= 0.122268182901557,\qquad & 
    a_3 &= 0.203023211433263, \qquad &
    b_1 &=  0.055200549768959, \\
    b_2 &= 0.127408150658963, &
    c_3 &= 0.001487834491987, & 
    \gamma_{\min} &\approx 0.0000189,\; \mathrm{Eff} \approx 21.98.
\end{alignedat}
\end{equation}

\noindent\textit{Variant 5 ($n_f = 7$, $p=4$, $v=3$).} Putting $c_1 = 0$ and $c_3 = 0$, one finds 
\begin{equation}\label{eq:BADABABADAB}
\begin{alignedat}{3}
    a_2 &= 0.201110227930330,\qquad & 
    a_3 &= 0.200577842713366, \qquad &
    b_1 &= 0.065692416344302, \\ 
    b_2 &= 0.264163604920340, &
    c_2 &= 0.001036943019757, & 
    \gamma_{\min} &\approx 0.00000520,\; \mathrm{Eff} \approx 80.13.
\end{alignedat}
\end{equation}

\noindent\textit{Variant 6 ($n_f = 6$, $p=4$, $v=3$).} For $c_2 = 0$ and $c_3 = 0$, minimization of $\mathrm{Err}_5$ results in
\begin{equation}\label{eq:DABABABABAD}
\begin{alignedat}{3}
    a_2 &= 0.282918304065611,\qquad & 
    a_3 &= -0.002348009438292, \qquad &
    b_1 &= 0.080181913812571, \\
    b_2 &= -1.372969015964262, &
    c_1 &= 0.000325098077953, & 
    \gamma_{\min} &\approx 0.0000166,\; \mathrm{Eff} \approx 46.47. 
\end{alignedat}
\end{equation}

\noindent\textit{Variant 7 ($n_f = 5$, $p=4$, $v=2$).} Finally, letting $c_1 = c_2 = c_3 = 0$, we obtain the optimal eleven-stage non-gradient scheme (in velocity form) by choosing
\begin{equation}\label{eq:BABABABABAB}
\begin{alignedat}{3}
    a_2 &= 0.253978510841060,\qquad & 
    a_3 &= -0.032302867652700, \qquad &
    b_1 &= 0.083983152628767, \\ 
    b_2 &=  0.682236533571909, &
    \gamma_{\min} &\approx 0.0000270, & \mathrm{Eff} &\approx 59.26. 
\end{alignedat}
\end{equation}

\subsubsection{Position versions}
When starting with a position update \eqref{eq:FGI_pos-update}, the eleven-stage Hessian-free FGI reads 
\begin{align*}
    \Psi_h &= \exp\left(a_1 h \hat{\T}\right) \exp\left(b_1 h \hat{\D}(h,b_1,c_1)\right) \exp\left(a_2 h \hat{\T}\right) \exp\left(b_2 h \hat{\D}(h,b_2,c_2)\right) \exp\left((\tfrac{1}{2} - (a_1+a_2))h \hat{\T}\right) \\
    &\qquad \times \exp\left((1-2(b_1+b_2))h \hat{\D}(h,1-2(b_1+b_2),c_3)\right) \exp\left((\tfrac{1}{2} - (a_1+a_2))h \hat{\T}\right) \\
    &\qquad \times \exp\left(b_2 h \hat{\D}(h,b_2,c_2)\right) \exp\left(a_2 h \hat{\T}\right) \exp\left(b_1 h \hat{\D}(h,b_1,c_1)\right) \exp\left(a_1 h \hat{\T}\right),
\end{align*}
where $a_1,a_2,b_1,b_2,c_1,c_2,c_3$ are free parameters.

\noindent\textit{Variant 1 ($n_f = 10$, $p=6$, $v=0$).} In general, there exist four real solutions to obtain a Hessian-free FGI of order $p=6$. The best solution is 
\begin{equation}\label{eq:ADADADADADA}
\begin{alignedat}{3}
    a_1 &= 0.109534125980058,\qquad & 
    a_2 &= 0.426279051773841, \qquad &
    b_1 &= 0.268835839917653, \\
    b_2 &= 0.529390037396794, &
    c_1 &= 0.000806354602850, & 
    c_2 &= 0.007662601517364, \\
    c_3 &=-0.011627206142396, & 
    \zeta_{\min} &\approx 0.00000603, & \mathrm{Eff} &\approx 0.17.
\end{alignedat}
\end{equation}
All four solutions have coefficients $c_j \neq 0$, i.e., an order reduction appears as soon as one $c_j$ is supposed to be zero. The variant with $c_1 = 0$ is omitted since order $p=6$ can be achieved by setting $a_1 = 0$, resulting in the nine-stage velocity algorithm \eqref{eq:BADADADAB}. 

\noindent\textit{Variant 2 ($n_f = 8$, $p=4$, $v=4$).} Putting $c_2 = 0$ yields
\begin{equation}\label{eq:ADABADABADA}
\begin{alignedat}{3}
    a_1 &= 0.083684971641549,\qquad & 
    a_2 &= 0.225966488946428, \qquad &
    b_1 &= 0.199022868372193, \\
    b_2 &= 0.197953981691206, &
    c_1 &= 0.000437056543403, & 
    c_3 &= 0.000870457820984, \\ 
    \gamma_{\min} &\approx 0.00000318, & \mathrm{Eff} &\approx 76.79. & &
\end{alignedat}
\end{equation}

\noindent\textit{Variant 3 ($n_f = 9$, $p=4$, $v=4$).} For $c_3 = 0$, the minimum is achieved at
\begin{equation}\label{eq:ADADABADADA}
\begin{alignedat}{3}
    a_1 &= 0.082541033171754,\qquad & 
    a_2 &= 0.228637847036999, \qquad &
    b_1 &= 0.196785139280847, \\
    b_2 &= 0.206783248777282, &
    c_1 &= 0.000317260402502, & 
    c_2 &= 0.000555360763892, \\ 
    \gamma_{\min} &\approx 0.00000235, & \mathrm{Eff} &\approx 64.99. & &
\end{alignedat}
\end{equation}

\noindent\textit{Variant 4 ($n_f = 6$, $p=4$, $v=3$).} Setting $c_1 = 0$ and $c_2 = 0$ results in
\begin{equation}\label{eq:ABABADABABA}
\begin{alignedat}{3}
    a_1 &= 0.134257092137626,\qquad & 
    a_2 &=-0.007010267216916, \qquad &
    b_1 &=-0.485681409840328, \\
    b_2 &= 0.767464037573892, &
    c_3 &= 0.002836723107629, & 
    \gamma_{\min} &\approx 0.0000154, \; \mathrm{Eff} \approx 50.09.
\end{alignedat}
\end{equation}

\noindent\textit{Variant 5 ($n_f = 7$, $p=4$, $v=3$).} At $c_1 = 0$ and $c_3 = 0$, the best solution reads
\begin{equation}\label{eq:ABADABADABA}
\begin{alignedat}{3}
    a_1 &= 0.062702644098210,\qquad & 
    a_2 &= 0.193174566017780, \qquad &
    b_1 &= 0.149293739165427, \\
    b_2 &= 0.220105234408407, &
    c_2 &= 0.000966194415594, & 
    \gamma_{\min} &\approx 0.00000445,\; \mathrm{Eff} \approx 93.60.
\end{alignedat}
\end{equation}

\noindent\textit{Variant 6 ($n_f = 7$, $p=4$, $v=3$).} Letting $c_2 = 0$ and $c_3 = 0$, the minimization of $\mathrm{Err}_5$ leads to
\begin{equation}\label{eq:ADABABABADA}
\begin{alignedat}{3}
    a_1 &= 0.115889910143319,\qquad & 
    a_2 &= 0.388722377182381,  \qquad &
    b_1 &= 0.282498420841510, \\
    b_2 &=-0.625616553474143, &
    c_1 &= 0.001208219887746, & 
    \gamma_{\min} &\approx 0.0000128,\; \mathrm{Eff} \approx 32.64.
\end{alignedat}
\end{equation}

\noindent\textit{Variant 7 ($n_f = 5$, $p=4$, $v=2$).} Demanding $c_1 = c_2 = c_3 = 0$, we arrive at the non-gradient scheme 
\begin{equation}\label{eq:ABABABABABA}
\begin{alignedat}{3}
    a_1 &= 0.275008121233242,\qquad & 
    a_2 &=-0.134795009910679, \qquad &
    b_1 &=-0.084429619507071, \\
    b_2 &= 0.354900057157426, &
    \gamma_{\min} &\approx 0.0000518, & \mathrm{Eff} &\approx 30.89.
\end{alignedat}
\end{equation}

\subsubsection{Remarks} The variants include two additional Hessian-free FGIs of order $p=6$, a velocity version \eqref{eq:BADADADADAB} with efficiency $\mathrm{Eff} \approx 0.27$ and a position version \eqref{eq:ADADADADADA} with $\mathrm{Eff} \approx 0.17$. Both versions are significantly more efficient than the nine-stage algorithm \eqref{eq:BADADADAB} with $\mathrm{Eff} \approx 0.0055$: using the same number of total force evaluations, the eleven-stage algorithms reduce the global errors approximately by a factor of $0.27/0.0055 \approx 49.09$ and $0.17/0.0055 \approx 30.91$, respectively. %Compared to the order-6 FGIs in \cite{omelyan2003symplectic}, the adapted FGIs \eqref{eq:BADADADADAB} and \eqref{eq:ADADADADADA} are more efficient.

The most efficient eleven-stage Hessian-free FGI of order $p=4$ is given by the position algorithm \eqref{eq:ABADABADABA} with $\mathrm{Eff} \approx 93.60$. Compared to the most efficient nine-stage algorithm \eqref{eq:BADABADAB} of order four ($\mathrm{Eff} \approx 73.45$), the eleven-stage algorithm reduces the global errors approximately by a factor of $93.60/73.45 \approx 1.27$. The most efficient velocity algorithm is \eqref{eq:BADABABADAB} with an efficiency of approximately $80.13$.
The non-gradient algorithm \eqref{eq:BABABABABAB} ($\mathrm{Eff} \approx 59.3$) became a state-of-the-art scheme in lattice quantum chromodynamics \cite{knechtli2017lattice} due to its high efficiency compared to lower-stage algorithms. Based on the efficiency measure, both the position algorithm \eqref{eq:ABADABADABA} and the velocity algorithm \eqref{eq:BABABABABAB} are able to reduce the global errors while keeping the number of total force recalculations fixed, namely by a factor of approximately $1.58$ and $1.35$, respectively.

\subsection{Complete classification}
A summary of all schemes derived in this section is provided in Table \ref{tab:aFGIs}. The table contains the schemes in abbreviated forms, i.e., $A$ and $B$ correspond to the exponential operators $\exp(a_j h \hat{\T})$ and $\exp(b_j h \hat{\V})$, respectively, and $D$ denotes approximated FG-updates of the form $\exp(b_j h \hat{\D}(h,b_j,c_j))$. For any number of stages $s=3,5,7,9,11$, the integrators are grouped according to their order of convergence. Within each group, the schemes are sorted in ascending order with respect to the number of force evaluations $n_f$. For the same number of evaluations, more efficient schemes are listed first. While some of the new schemes have only slightly different time coefficients compared to the results in \cite{omelyan2003symplectic}, there are also variants where one can observe large differences in the time coefficients.

\begin{table}
    \centering
    \begin{tabular}{l c r r r r c r}
        \toprule
        Algorithm & Equation & $p$ & $n_f$ & $\mathrm{Err}_{p+1}$ & Efficiency & Remarks & ID in \cite{omelyan2003symplectic} \\
        \midrule
        BAB & \eqref{eq:BAB} & 2 & 1 & 0.0932 & 10.73 & \cite{verlet1967computer} & 1\\
        ABA & \eqref{eq:ABA} & 2 & 1 & 0.0932 & 10.73 & \cite{verlet1967computer} & 2\\
        DAD & \eqref{eq:DAD} & 2 & 2 & 0.0833 & 3.00 & \cite{omelyan2003symplectic,Hairer_McLachlan_Skeel_2009} & 3\\
        ADA & \eqref{eq:ADA} & 2 & 2 & 0.0417 & 6.00 & \cite{omelyan2003symplectic} & 4\\
        \midrule
        BABAB & \eqref{eq:BABAB} & 2 & 2 & 0.00855 & 29.24 & \cite{mclachlan1995,omelyan2003symplectic} & 5\\
        ABABA & \eqref{eq:ABABA} & 2 & 2 & 0.00855 & 29.24 & \cite{mclachlan1995,omelyan2003symplectic} & 6\\
        BADAB & \eqref{eq:BADAB} & 4 & 3 & 0.000728 &16.96 & \cite{suzuki1995new,chin1997symplectic} & 8\\
        DABAD & \eqref{eq:DABAD} & 4 & 3 & 0.00335 &  3.68 & \cite{omelyan2003symplectic} & 7\\
        DADAD & \eqref{eq:DADAD} & 4 & 4 & 0.000625 & 6.25 & New & 9\\
        ADADA & \eqref{eq:ADADA} & 4 & 4 & 0.000718 & 5.44 & \cite{suzuki1995new,chin1997symplectic} & 10\\
        \midrule
        ABABABA & \eqref{eq:ABABABA} & 4 & 3 & 0.0283 & 0.44 & \cite{forest1990fourth,yoshida1990construction} & 12\\
        BABABAB & \eqref{eq:BABABAB} & 4 & 3 & 0.0383 & 0.32 & \cite{forest1990fourth,yoshida1990construction} & 11\\
        ABADABA & \eqref{eq:ABADABA} & 4 & 4 & 0.000149 & 26.19 & New & 14\\
        DABABAD & \eqref{eq:DABABAD} & 4 & 4 & 0.000891 & 4.38 & New & 13\\
        BADADAB & \eqref{eq:BADADAB} & 4 & 5 & 0.0000498 & 32.12 & New & 15\\
        ADABADA & \eqref{eq:ADABADA} & 4 & 5 & 0.0000844 & 18.95 & New & 16\\
        ADADADA & \eqref{eq:ADADADA} & 4 & 6 & 0.0000200 & 38.57 & New & 18\\
        DADADAD & \eqref{eq:DADADAD} & 4 & 6 & 0.0000275 & 28.09 & New & 17\\
        \midrule
        ABABABABA & \eqref{eq:ABABABABA} & 4 & 4 & 0.000610 & 6.40 & \cite{omelyan2003symplectic} & 20\\
        BABABABAB & \eqref{eq:BABABABAB} & 4 & 4 & 0.000654 & 5.97 & \cite{omelyan2003symplectic} & 19\\
        BABADABAB & \eqref{eq:BABADABAB} & 4 & 5 & 0.0000651 & 24.57 & New & 21\\
        DABABABAD & \eqref{eq:DABABABAD} & 4 & 5 & 0.000336 & 4.76 & New & 22\\
        BADABADAB & \eqref{eq:BADABADAB} & 4 & 6 & 0.0000105 & 73.45 & New & 24\\
        DABADABAD & \eqref{eq:DABADABAD} & 4 & 6 & 0.0000130 & 59.33 & New & 23\\
        ABADADABA & \eqref{eq:ABADADABA} & 4 & 6 & 0.0000346 & 22.32 & New & 25\\
        ADABABADA & \eqref{eq:ADABABADA} & 4 & 6 & 0.0000471 & 16.39 & New & 26\\
        DADABADAD & \eqref{eq:DADABADAD} & 4 & 7 & 0.0000101 & 41.06 & New & 27\\
        ADADADADA & \eqref{eq:ADADADADA} & 4 & 8 & 0.00000501 & 48.71 & New & 29\\
        BADADADAB & \eqref{eq:BADADADAB} & 6 & 7 & 0.00154 & 0.0055 & \cite{omelyan2003symplectic} & 28\\
        \midrule
        BABABABABAB & \eqref{eq:BABABABABAB} & 4 & 5 & 0.0000270 & 59.26 & \cite{omelyan2003symplectic} & 30\\
        ABABABABABA & \eqref{eq:ABABABABABA} & 4 & 5 & 0.0000518 & 30.89 & \cite{omelyan2003symplectic} & 31\\
        ABABADABABA & \eqref{eq:ABABADABABA} & 4 & 6 & 0.0000154 & 50.09 & New & 33\\
        DABABABABAD & \eqref{eq:DABABABABAD} & 4 & 6 & 0.0000166 & 46.47 & New & 32\\
        ABADABADABA & \eqref{eq:ABADABADABA} & 4 & 7 & 0.00000445 & 93.60 & New & 36\\
        BADABABADAB & \eqref{eq:BADABABADAB} & 4 & 7 & 0.00000520 & 80.13 & New & 34\\
        ADABABABADA & \eqref{eq:ADABABABADA} & 4 & 7 & 0.0000128 & 32.64 & New & 37\\
        BABADADABAB & \eqref{eq:BABADADABAB} & 4 & 7 & 0.0000189 & 21.98 & New & 35\\
        ADABADABADA & \eqref{eq:ADABADABADA} & 4 & 8 & 0.00000318 & 76.79 & New & 40\\
        DABADADABAD & \eqref{eq:DABADADABAD} & 4 & 8 & 0.00000355 & 68.84 & New & 38\\
        DADABABADAD & \eqref{eq:DADABABADAD} & 4 & 8 & 0.00000519 & 47.08 & New & 39\\
        ADADABADADA & \eqref{eq:ADADABADADA} & 4 & 9 & 0.00000235 & 64.99 & New & 43\\
        BADADADADAB & \eqref{eq:BADADADADAB} & 6 & 9 & 0.00000699 & 0.27 & New & 42\\
        ADADADADADA & \eqref{eq:ADADADADADA} & 6 & 10 & 0.00000603 & 0.17 & New & 45\\
        \bottomrule
    \end{tabular}
    \caption{Collection of Hessian-free force-gradient integrators with up to eleven stages. For all algorithms, the convergence order $p$, the number of force evaluations $n_f$ per time step, the leading error term $\mathrm{Err}_{p+1}$, its efficiency $\mathrm{Eff}$, as well as a remark on the earlier appearance of the time coefficients, are stated. For completeness, the final column refers to the ID of the respective variant in \cite{omelyan2003symplectic}.}
    \label{tab:aFGIs}
\end{table}

\subsection{Hessian-free force-gradient integrators of higher orders}
In case of very accurate integrations, it can be useful to consider numerical time integration schemes of even higher convergence order $p$. The derivation of Hessian-free and exact FGIs via direct decompositions will result in the most efficient schemes within the framework. However, the derivation becomes a more challenging numerical problem when further increasing the number of stages. 

An alternative approach of deriving numerical time integration schemes of arbitrarily high convergence order is given by means of composition techniques \cite{yoshida1990construction,suzuki1990fractal,omelyan2002construction}. Using symmetric compositions, the time-reversibility, volume-preservation and closure property of the underlying base methods is preserved.

Since the derivation based on direct decomposition will always lead to more efficient algorithms, composition techniques should only be applied for very high orders for which no methods based on direct decomposition are known. For example, Yoshida's triple jump \cite{yoshida1990construction} based on the non-gradient scheme \eqref{eq:BAB} yields the method \eqref{eq:BABABAB} which is the least effective decomposition scheme ($\mathrm{Eff} \approx 0.32$) of order four under investigation. 

Starting from a base method of order $p$, Yoshida's triple jump requires $3^k$ applications of the base method to obtain a composition scheme of order $p+2k$. The alternative approach of Suzuki's fractals \cite{suzuki1990fractal} demands even more applications, namely $5^k$. Advanced composition schemes \cite{omelyan2002construction} provide a possibility to derive composition schemes with significantly less stages if one wants to increase the order by more than two orders. Advanced composition schemes starting from a base scheme of order $p \in \{4,6,8\}$ has been discussed in \cite{omelyan2002construction}. Advanced composition schemes based on base schemes of order $p=2$ have been investigated in \cite{kahan1997composition}.

\section{Numerical results}\label{sec:Numerical_Results}
This section provides numerical tests of Hessian-free FGIs for three different test examples. In Section \ref{sec:OSS}, numerical tests for the outer solar system, a $N$-body problem with mass matrix $\bfM \neq \mathrm{Id}$, are performed. Then, in Section \ref{sec:HMC}, we consider the Hybrid Monte Carlo (HMC) algorithm \cite{duane1987hybrid} in two different applications. The two-dimensional Schwinger model is discussed in Section \ref{sec:Schwinger}. Here, we have implemented the analytical FG-term \cite{shcherbakov2017adapted}. This enables a comparison of the exact FGIs \cite{omelyan2003symplectic} to the new class of Hessian-free FGIs. Secondly, in Section \ref{sec:QCD}, four-dimensional gauge field simulations in lattice QCD with two heavy Wilson fermions are discussed. Since the links are elements of the special unitary group $\mathrm{SU}(3)$, this provides an example in a non-Abelian setting.

\subsection{The outer solar system}\label{sec:OSS}
We consider the separable Hamiltonian system 
\begin{equation*}%\label{eq:OSS}
    \Ham(\p,\q) =  \frac{1}{2}\sum\limits_{i=0}^5 \frac{1}{m_i} p_i^\top p_i - G \sum\limits_{i=1}^5 \sum\limits_{j=0}^{i-1} \frac{m_i m_j}{\lVert q_i - q_j \rVert},
\end{equation*}
describing the motion of the four outer planets (and Pluto) relative to the sun. Here, the momenta $\p$ and positions $\q$ are supervectors composed by the vectors $p_i,q_i \in \mathbb{R}^3$. We have taken the parameters from \cite{HairerLubichWanner} and also took the initial data from \cite{burckhardt1993ahnerts} corresponding to September 5, 1994 at midnight. 
Solutions of the outer solar system using the Hessian-free FGI \eqref{eq:ABADABADABA} with step size $h=200$ over a time period of 200,000 days are depicted in Figure \ref{fig:OSS_sol}, showing the correct behavior. Numerical results for the other variants of Hessian-free FGIs show the same behavior, as expected.

\begin{minipage}[t]{0.39\textwidth}
\begin{figure}[H]
    \centering
    \vspace{-22.25pt}
    \scalebox{0.8}{\input{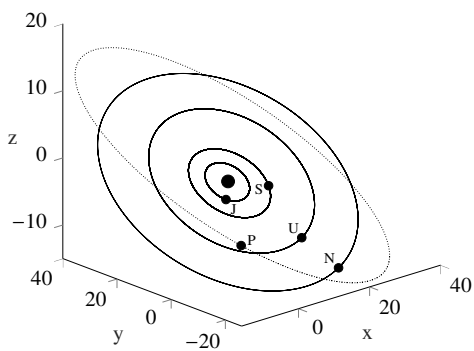}}
    \vspace{-12pt}
    \caption{Outer solar system. Numerical solution using the Hessian-free FGI \eqref{eq:ABADABADABA} and $h=200$ over a time period of 200,000 days.}
    \label{fig:OSS_sol}
\end{figure}
\end{minipage} \hfill 
\begin{minipage}[t]{0.55\textwidth}
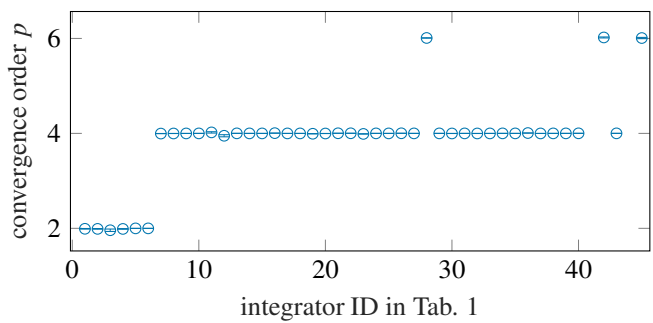
\begin{figure}[H]
    \centering
    % This file was created by matlab2tikz.
%
%The latest updates can be retrieved from
%  http://www.mathworks.com/matlabcentral/fileexchange/22022-matlab2tikz-matlab2tikz
%where you can also make suggestions and rate matlab2tikz.
%
\definecolor{mycolor1}{rgb}{0.00000,0.46275,0.68627}%
\begin{tikzpicture}

\begin{axis}[%
width=3.0in,
height=1.25in,
at={(0.4in,0.423in)},
scale only axis,
xmin=-0.144031512182604,
xmax=45.6369629875944,
xlabel style={font=\color{white!15!black}},
xlabel={integrator ID in Tab. \ref{tab:aFGIs}},
ymin=1.5235648623295,
ymax=6.56709911152467,
ylabel style={font=\color{white!15!black}},
ylabel={convergence order $p$},
ylabel shift = -1pt,
axis background/.style={fill=white},
legend style={legend cell align=left, align=left, draw=white!15!black},
every axis legend/.code={\let\addlegendentry\relax}
]
\addplot [color=mycolor1, draw=none, mark=o, mark options={solid, mycolor1}]
 plot [error bars/.cd, y dir = both, y explicit]
 table[row sep=crcr, y error plus index=2, y error minus index=3]{%
1	1.99035443434003	0.00610066711489857	0.00610066711489857\\
2	1.98867260055869	0.00693843091934269	0.00693843091934269\\
3	1.95966126823971	0.025533913989938	0.025533913989938\\
4	1.98753433616851	0.00706895840196012	0.00706895840196012\\
5	1.99887965367715	0.000624395775704249	0.000624395775704249\\
6	1.99908531401874	0.000555946736211637	0.000555946736211637\\
7	3.99594273358794	0.00174134298039135	0.00174134298039135\\
8	3.99787804503562	0.000980226369449536	0.000980226369449536\\
9	3.99884908704969	0.000562060478363348	0.000562060478363348\\
10	3.99959019652774	0.000214671545084318	0.000214671545084318\\
11	4.02077545451398	0.0140714526131093	0.0140714526131093\\
12	3.94992153520423	0.0222139883375344	0.0222139883375344\\
13	4.0005419650175	0.000252825436956343	0.000252825436956343\\
14	3.999564211525	0.000218282443548339	0.000218282443548339\\
15	3.99929463398883	0.000338923884542677	0.000338923884542677\\
16	4.0050060723689	0.00209514998164089	0.00209514998164089\\
17	3.99926361167124	0.000353406031516105	0.000353406031516105\\
18	3.99930075325723	0.000345523166643367	0.000345523166643367\\
19	3.99215460009554	0.0035008338302113	0.0035008338302113\\
20	3.99745341076421	0.00116923994636163	0.00116923994636163\\
21	4.00222327210791	0.000890413898798689	0.000890413898798689\\
22	4.00196720554842	0.000873574126318449	0.000873574126318449\\
23	3.98873652806332	0.00446115923085448	0.00446115923085448\\
24	4.0000114056025	7.45516801688687e-05	7.45516801688687e-05\\
25	3.99954419442251	0.000250170651273975	0.000250170651273975\\
26	4.00194768239466	0.000734987048812044	0.000734987048812044\\
27	3.99969946912994	0.000200464231145319	0.000200464231145319\\
28	6.00971510643407	0.00422834599499522	0.00422834599499522\\
29	3.99990991970191	0.000165620198875671	0.000165620198875671\\
30	3.99819143344861	0.000813300573932933	0.000813300573932933\\
31	3.99962384910811	0.000198499197910882	0.000198499197910882\\
32	3.99837706847127	0.000734178617243936	0.000734178617243936\\
33	3.99918530677084	0.000390245075319871	0.000390245075319871\\
34	4.00009333319857	9.82750846402569e-05	9.82750846402569e-05\\
35	3.99970297994007	0.000208332077894178	0.000208332077894178\\
36	4.0071403793787	0.00210279566422596	0.00210279566422596\\
37	3.99881270516463	0.000550621257231085	0.000550621257231085\\
38	3.99971409279748	0.00020371819738063	0.00020371819738063\\
39	4.00002181435353	0.00013114017005608	0.00013114017005608\\
40	4.00114628387279	0.000185846459588825	0.000185846459588825\\
41	0	0	0\\
42	6.02045611202353	0.010209476085981	0.010209476085981\\
43	4.00054757392147	0.000197187006086842	0.000197187006086842\\
44	0	0	0\\
45	6.00822514945997	0.0133018070299698	0.0133018070299698\\
};
\addlegendentry{data1}

\end{axis}

\begin{axis}[%
width=3.5in,
height=1.75in,
at={(0in,0in)},
scale only axis,
xmin=0,
xmax=1,
ymin=0,
ymax=1,
axis line style={draw=none},
ticks=none,
axis x line*=bottom,
axis y line*=left,
legend style={legend cell align=left, align=left, draw=white!15!black}
]
\end{axis}
\end{tikzpicture}%
    \caption{Outer solar system. Numerical verification of the convergence order of Hessian-free FGIs. For each variant of Hessian-free FGIs, the scaling has been measured performing simulations using different step sizes. The plot displays the mean, as well as the standard deviation of the numerical measurement of the convergence order.}
    \label{fig:OSS_convergence-order}
\end{figure}
\end{minipage}
\vspace*{12pt}

Figure \ref{fig:OSS_convergence-order} shows that all Hessian-free FGIs achieve the desired convergence order. This results are more or less just a reference that the derivation and solution of the order conditions has been performed correctly and that Hessian-free FGIs are also applicable to Hamiltonian systems \eqref{eq:Hamiltonian_general} with $\bfM \neq \mathrm{Id}$.

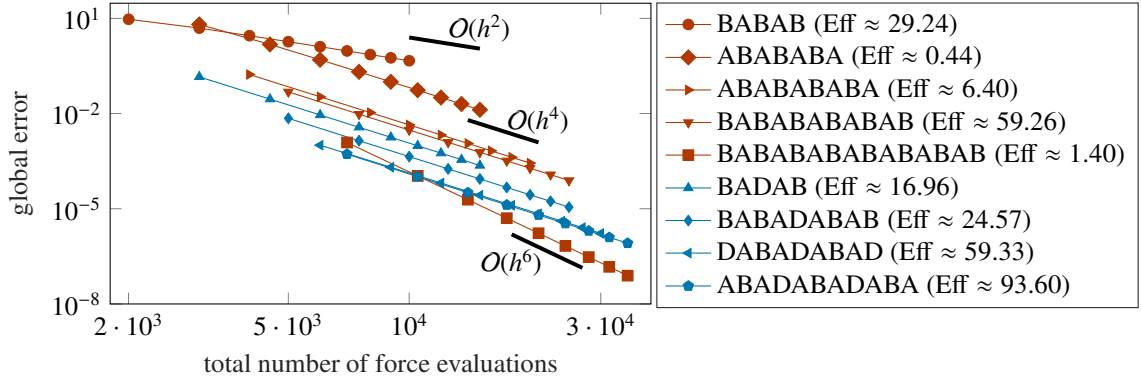
\begin{figure}[htb]
    \centering
    % This file was created by matlab2tikz.
%
%The latest updates can be retrieved from
%  http://www.mathworks.com/matlabcentral/fileexchange/22022-matlab2tikz-matlab2tikz
%where you can also make suggestions and rate matlab2tikz.
%
\definecolor{mycolor1}{rgb}{0.68627,0.22353,0.00000}%
\definecolor{mycolor2}{rgb}{0.00000,0.46275,0.68627}%
\begin{tikzpicture}

\begin{axis}[%
width=2.8in,
height=1.57in,
at={(0.575in,0.45in)},
scale only axis,
xmode=log,
xmin=1800,
xmax=40000,
xminorticks=true,
xtick={2000,5000,10000,30000},
xticklabels={$2\cdot 10^3$,$5\cdot 10^3$,$10^4$,$3\cdot 10^4$},
xlabel style={font=\color{white!15!black}},
xlabel={total number of force evaluations},
ymode=log,
ymin=1e-08,
ymax=30,
yminorticks=true,
ylabel style={font=\color{white!15!black}},
ylabel={global error},
axis background/.style={fill=white},
legend style={at={(1.01,1.0025)}, anchor=north west, legend cell align=left, align=left, draw=white!15!black}
]
\addplot [color=black, line width=1.5pt, forget plot]
  table[row sep=crcr]{%
14000	0.00625\\
17500	0.00256\\
21000	0.00123456790123457\\
};
\node[] at (axis cs: 2.1e+4,6.25e-3) {$\mathcal{O}(h^4)$};

\addplot [color=black, line width=1.5pt, forget plot]
  table[row sep=crcr]{%
10000	2.5\\
12500	1.6\\
15000	1.11111111111111\\
};
\node[] at (axis cs: 1.5e+4,0.5e+1) {$\mathcal{O}(h^2)$};

\addplot [color=black, line width=1.5pt, forget plot]
  table[row sep=crcr]{%
18000	1.5625e-06\\
22500	4.096e-07\\
27000	1.37174211248285e-07\\
};
\node[] at (axis cs: 1.8e+4,2e-7) {$\mathcal{O}(h^6)$};

\addplot [color=mycolor1,   mark=*, mark options={solid, mycolor1}]
  table[row sep=crcr]{%
2001	9.36660756694869\\
3001	4.92916439680172\\
4001	2.8518766239977\\
5001	1.83950330055993\\
6001	1.28111479698173\\
7001	0.942427731591893\\
8001	0.722009823817705\\
9001	0.570680204769919\\
10001	0.462349423708933\\
};
\addlegendentry{BABAB ($\mathrm{Eff} \approx 29.24$)}

\addplot [color=mycolor1,   mark=square*,  mark options={solid, mycolor1,rotate=45}]
  table[row sep=crcr]{%
3000	6.46388628774621\\
4500	1.50986048556565\\
6000	0.495576364738267\\
7500	0.206290915164151\\
9000	0.100351119234063\\
10500	0.0544500101195447\\
12000	0.0320255160707186\\
13500	0.0200396417636508\\
15000	0.0131697495266549\\
};
\addlegendentry{ABABABA ($\mathrm{Eff} \approx 0.44$)}

\addplot [color=mycolor1,   mark=triangle*, mark options={solid, mycolor1,rotate=270}]
  table[row sep=crcr]{%
4000	0.170186597952258\\
6000	0.0338239969585849\\
8000	0.0107223692249845\\
10000	0.00439552894064838\\
12000	0.00212068767548239\\
14000	0.0011449923154448\\
16000	0.000671285413831277\\
18000	0.000419127497093063\\
20000	0.000275011489117369\\
};
\addlegendentry{ABABABABA ($\mathrm{Eff} \approx 6.40$)}

\addplot [color=mycolor1,   mark=triangle*, mark options={solid, rotate=180, mycolor1}]
  table[row sep=crcr]{%
5001	0.0482826817171351\\
7501	0.00957433227144804\\
10001	0.00303328939418458\\
12501	0.00124316259261889\\
15001	0.000599707260837482\\
17501	0.000323767669362428\\
20001	0.000189809330021817\\
22501	0.000118506430991012\\
25001	7.77562857943315e-05\\
};
\addlegendentry{BABABABABAB ($\mathrm{Eff} \approx 59.26$)}

\addplot [color=mycolor1, mark=square*, mark options={solid, mycolor1}]
  table[row sep=crcr]{%
7001	0.00122378410215485\\
10501	0.000108148733804426\\
14001	1.92901030676358e-05\\
17501	5.06098817434012e-06\\
21001	1.69500947387729e-06\\
24501	6.71632566812694e-07\\
28001	3.00803150195671e-07\\
31501	1.47730780392454e-07\\
35001	7.79621074811832e-08\\
};
\addlegendentry{BABABABABABABAB ($\mathrm{Eff} \approx 1.40$)}

\addplot [color=mycolor2, mark=triangle*, mark options={solid, mycolor2}]
  table[row sep=crcr]{%
3001	0.143684589860833\\
4501	0.0285303594689209\\
6001	0.00904151368097477\\
7501	0.00370597102268692\\
9001	0.00178786888860871\\
10501	0.000965255708863546\\
12001	0.000565893394686194\\
13501	0.000353317441415464\\
15001	0.000231826773967143\\
};
\addlegendentry{BADAB ($\mathrm{Eff} \approx 16.96$)}

\addplot [color=mycolor2, mark=diamond*, mark options={solid, mycolor2}]
  table[row sep=crcr]{%
5001	0.00709588660183641\\
7501	0.00139580938350111\\
10001	0.00044098617185175\\
12501	0.000180502248480099\\
15001	8.70143188945666e-05\\
17501	4.69568796778624e-05\\
20001	2.75206720747664e-05\\
22501	1.71788098517773e-05\\
25001	1.12697629288784e-05\\
};
\addlegendentry{BABADABAB ($\mathrm{Eff} \approx 24.57$)}

\addplot [color=mycolor2, mark=triangle*, mark options={solid, rotate=90, mycolor2}]
  table[row sep=crcr]{%
6002	0.00100779141940163\\
9002	0.000203466462711674\\
12002	6.48679365394875e-05\\
15002	2.66634922681549e-05\\
18002	1.28835546739888e-05\\
21002	6.96271403331008e-06\\
24002	4.08491757857151e-06\\
27002	2.55197682812526e-06\\
30002	1.67535098986684e-06\\
};
\addlegendentry{DABADABAD ($\mathrm{Eff} \approx 59.33$)}

\addplot [color=mycolor2, mark=pentagon*, mark options={solid, mycolor2}]
  table[row sep=crcr]{%
7000	0.000525371143029686\\
10500	0.000102510154156832\\
14000	3.22933531377223e-05\\
17500	1.319993642521e-05\\
21000	6.35805016293897e-06\\
24500	3.42906159562379e-06\\
28000	2.00865724077753e-06\\
31500	1.25317322704096e-06\\
35000	8.21520482886588e-07\\
};
\addlegendentry{ABADABADABA ($\mathrm{Eff} \approx 93.60$)}

\end{axis}

\begin{axis}[%
width=6in,
height=2.1in,
at={(0in,0in)},
scale only axis,
xmin=0,
xmax=1,
ymin=0,
ymax=1,
axis line style={draw=none},
ticks=none,
axis x line*=bottom,
axis y line*=left,
legend style={legend cell align=left, align=left, draw=white!15!black}
]
\end{axis}
\end{tikzpicture}%
    \caption{Outer solar system. Global error vs.\ number of force evaluations for a selected number of Hessian-free FGIs (blue lines) and non-gradient schemes (red lines). The simulation has been performed over a time period of 200,000 days using different number of steps, $\mathrm{nsteps} \in \{1000,1500,\ldots,5000\}$.}
    \label{fig:OSS_efficiency}
\end{figure}

Furthermore, the efficiency of the integrators for the outer solar system is depicted in Figure \ref{fig:OSS_efficiency} by plotting the global error vs. the total number of force evaluations for different step sizes. 
For the sake of clarity, only a selection of integrators is shown: the most efficient non-gradient scheme for any number of stages, the frequently used Hessian-free FGI \eqref{eq:BADAB}, as well as the three best performing Hessian-free FGIs of order four.
Compared to the most efficient fourth-order non-gradient scheme \eqref{eq:BABABABABAB}, Hessian-free FGIs are able to reduce the error by more than one order of magnitude at the same number of force evaluations, emphasizing the efficiency of the proposed class of integrators.
For accuracy demands up to $1e-4$, the Hessian-free FGI \eqref{eq:ABADABADABA} is even able to compete with the most efficient sixth-order integrator under investigation, the 15-stage non-gradient velocity algorithm introduced in \cite{omelyan2003symplectic}. 

In addition, the results show that the practical efficiency of the integrators might diverge from the theoretical findings summarized in Table \ref{tab:aFGIs}. For example, \eqref{eq:BADAB} ($\mathrm{Eff} \approx 16.96$) is more efficient than \eqref{eq:BABABABABAB} ($\mathrm{Eff} \approx 59.26$) and \eqref{eq:DABADABAD} ($\mathrm{Eff} \approx 59.33$) is approximately as efficient as \eqref{eq:ABADABADABA} ($\mathrm{Eff} \approx 93.60$).
The theoretical efficiency values are derived by equating all brackets to one; therefore, this outcome is not contradictory. Nevertheless, the efficiency measure provides a valuable heuristic, indicating which integrators are likely to be efficient without requiring additional assumptions about the specific system.

\begin{figure}[htb]
    \centering
    \input{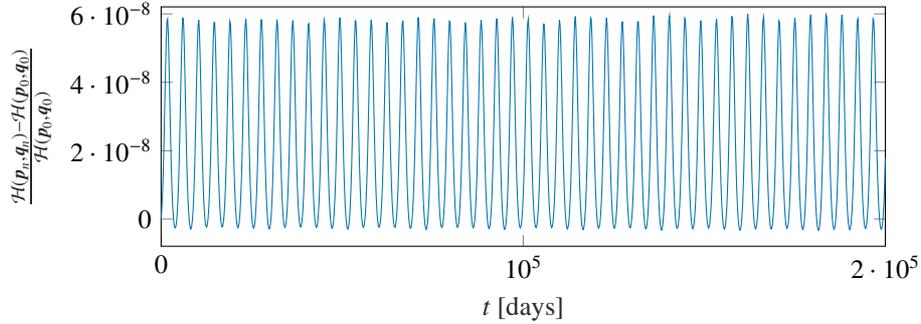}
    \caption{Outer solar system. Relative energy error of the Hessian-free FGI \eqref{eq:ABADABADABA} using $h=200$ over a time period of 200,000 days.}
    \label{fig:OSS_energy_error}
\end{figure}

Another important feature is the long-time energy conservation of Hessian-free FGIs. As an example, Figure \ref{fig:OSS_energy_error} depicts the relative energy error $(\Ham(\p_n,\q_n) - \Ham(\p_0,\q_0))/\Ham(\p_0,\q_0)$ for the Hessian-free FGI \eqref{eq:ABADABADABA}. For this problem, no drift is visible, i.e.\ the Hessian-free FGIs behave like a symplectic integrator and thus show excellent long-term energy behavior. Note that this observation coincides with the observations made in \cite{Hairer_McLachlan_Skeel_2009} for another $N$-body problem.

\subsection{Hybrid Monte Carlo simulations}\label{sec:HMC}
The most common approach of simulating quantum field theories on the lattice is given by Hybrid Monte Carlo (HMC) simulations \cite{duane1987hybrid}. Within the HMC algorithm, the lattice quantum field theory is embedded in a higher-dimensional classical system by introducing a fictitious simulation time. It is described by a separable Hamiltonian system
$$\Ham([\bfP],[\bfQ]) = \T([\bfP]) + \V([\bfQ]),$$
reaching the equilibrium distribution of the desired expectation values in the large time limit. The square brackets denote a configuration, i.e.\ the set of all variables arranged on the lattice. For a given configuration $[\bfQ]$, a new configuration $[\bfQ']$ is generated by performing an HMC update consisting of a molecular dynamics (MD) step and a Metropolis step.
\begin{itemize}
    \item \textbf{MD step.} For a trajectory of length $\tau$, solve the Hamiltonian equations of motion given by
    \begin{equation*}%\label{eq:md-step}
        \dot{\bfP}_{x,\mu} = - \frac{\partial \Ham([\bfP],[\bfQ])}{\partial \bfQ_{x,\mu}}, \qquad \dot{\bfQ}_{x,\mu} = \frac{\partial \Ham([\bfP],[\bfQ])}{\partial \bfP_{x,\mu}}
    \end{equation*}
    for each position $(x,\mu)$ on the lattice. We denote the new configuration by $([\bfP'],[\bfQ'])$.  
    \item \textbf{Metropolis step.} Accept the new configuration $([\bfP'],[\bfQ'])$ with probability 
    $$ \mathcal{P}([\bfQ] \to [\bfQ']) = \min\left(1 , \exp(- \Delta \Ham) \right),$$
    where $\Delta \Ham = \Ham([\bfP'],[\bfQ']) - \Ham([\bfP],[\bfQ])$.
    Depending on the acceptance step, either the new or old configuration is added to the Markov chain and used as the initial configuration for the next MD step.
\end{itemize}
The Markov process will converge to a unique fixed point distribution provided that it is ergodic and satisfies the detailed balance condition \cite{duane1987hybrid}. The detailed balance condition is satisfied, if the numerical integration scheme used in the MD step is time-reversible and volume-preserving \cite{knechtli2017lattice}. 
Good energy conservation is required to ensure a high acceptance rate in the Metropolis step. At the same time, one is interested in minimizing the computational cost of the HMC algorithm. Consequently, the choice of the integrator in the MD step plays a crucial role. The optimal acceptance rate for a numerical integration scheme of order $p$ is \cite{TAKAISHI20006} 
\begin{align}\label{eq:optimal_acceptance_rate}
    \langle P_{\mathrm{acc}}\rangle_{\mathrm{opt}} = \exp(-1/p).
\end{align}
Equation \eqref{eq:optimal_acceptance_rate} indicates that the optimal acceptance rate for order $p=2$ is $61\%$, $78\%$ for order $p=4$ and $85\%$ for order $p=6$. In practice, however, achieving the optimal acceptance rate is often challenging, as instabilities emerge as a limiting factor in the step size $h$. Hence one frequently aims at achieving an acceptance rate of $P_\mathrm{acc} \approx 90\%$.
For log-normal distributed $\Delta H$, the expected acceptance rate is given by \cite{GUPTA1990437}
\begin{align}\label{eq:acceptance_rate}
    \langle P_{\mathrm{acc}}\rangle = \mathrm{erfc} \left(\sqrt{\sigma^2/8}\right),
\end{align}
where $\sigma^2$ denotes the (empirical) variance in $\Delta H$ and $\mathrm{erfc}$ the complementary error function. All in all, we are interested in minimizing the number of force evaluations per unit trajectory $n_f / h$ while achieving a desired acceptance rate $P_\mathrm{acc} \in [0,1]$, i.e., we aim at maximizing the step size $h$ such that $\sigma^2 \leq 8 (\mathrm{erfc}^{-1}(P_\mathrm{acc}))^2$ is satisfied.

\subsubsection{The two-dimensional Schwinger model}\label{sec:Schwinger}
\begin{figure}[H]
    \centering
    \begin{subfigure}[b]{\textwidth}
        \centering
        % This file was created by matlab2tikz.
%
%The latest updates can be retrieved from
%  http://www.mathworks.com/matlabcentral/fileexchange/22022-matlab2tikz-matlab2tikz
%where you can also make suggestions and rate matlab2tikz.
%
\definecolor{mycolor1}{rgb}{0.68627,0.22353,0.00000}%
\definecolor{mycolor2}{rgb}{0.00000,0.46275,0.68627}%
\begin{tikzpicture}

\begin{axis}[%
width=5in,
height=0.85in,
at={(0.5in,0.425in)},
scale only axis,
xmin=0,
xmax=46,
ymin=0,
ymax=40,
xlabel style={font=\color{white!15!black}},
xlabel={integrator ID in \cite{omelyan2003symplectic}}, 
ylabel style={font=\color{white!15!black}},
ylabel={$n_f/h(P_{\mathrm{acc}} = 90\%)$},
axis background/.style={fill=white},
axis x line*=bottom,
axis y line*=left,
legend style={legend cell align=left, align=left, draw=white!15!black}
]
\addplot [color=mycolor1, draw=none, mark=o, mark options={solid, mycolor1}, forget plot]
 plot [error bars/.cd, y dir = both, y explicit]
 table[row sep=crcr, y error plus index=2, y error minus index=3]{%
1	14.9392169078503	0.941945657022855	0.941945657022855\\
3	36.80723294784	2.35744873173032	2.35744873173032\\
2	13.9756371363434	0.873248418847592	0.873248418847592\\
4	25.7382125900335	1.5157805885332	1.5157805885332\\
5	8.5238130144369	0.495827630804688	0.495827630804688\\
7	16.8807105027147	0.439780792430016	0.439780792430016\\
8	11.9388702595316	0.342900863452771	0.342900863452771\\
9	16.4334056016841	0.456775384564013	0.456775384564013\\
6	7.60447361415898	0.418355144116458	0.418355144116458\\
10	17.543055103488	0.470745354477339	0.470745354477339\\
11	23.2790961505983	0.673261684668855	0.673261684668855\\
13	21.5301741435329	0.583944380783353	0.583944380783353\\
15	8.7375198437765	0.234941607098367	0.234941607098367\\
17	9.96164847240225	0.272894348310837	0.272894348310837\\
12	22.1518630880617	0.664882937774305	0.664882937774305\\
14	9.209610526788	0.251927305725023	0.251927305725023\\
16	11.5555917950751	0.339024256708156	0.339024256708156\\
18	9.04434025247923	0.280939284314634	0.280939284314634\\
19	13.7405795668681	0.374068097701979	0.374068097701979\\
21	9.20463030288959	0.262152544246045	0.262152544246045\\
22	13.2452224335611	0.37539648560981	0.37539648560981\\
23	6.35718673505552	0.177729362578142	0.177729362578142\\
24	6.5628139369886	0.180501727849751	0.180501727849751\\
27	8.79567343581552	0.239058021218628	0.239058021218628\\
28	37.4361393552764	0.756102601949286	0.756102601949286\\
20	13.5673365844331	0.370817430459757	0.370817430459757\\
25	9.78276242467253	0.278707981151371	0.278707981151371\\
26	12.0708388808342	0.338197903730322	0.338197903730322\\
29	9.04631788375202	0.259344840112093	0.259344840112093\\
30	7.05768763542669	0.188392113126634	0.188392113126634\\
32	8.47386951220192	0.233052178173914	0.233052178173914\\
34	6.54874987351493	0.196446436361089	0.196446436361089\\
35	10.1104655213643	0.277302089324374	0.277302089324374\\
38	20.020431108263	0.375968390447116	0.375968390447116\\
39	25.2760985854508	0.469310619377317	0.469310619377317\\
42	24.4809473394992	0.494063282643668	0.494063282643668\\
44	25.6059106694543	0.536877278336468	0.536877278336468\\
31	9.18162134570112	0.255780104795702	0.255780104795702\\
33	8.45002350877251	0.252290328853316	0.252290328853316\\
36	5.28867375415944	0.150742315200507	0.150742315200507\\
37	9.65664439122953	0.28899378579775	0.28899378579775\\
40	20.6976743107926	0.394662753932568	0.394662753932568\\
41	29.9053523963048	0.625597849725263	0.625597849725263\\
43	27.8645043817555	0.524336959754409	0.524336959754409\\
45	27.0112776726202	0.492422524159978	0.492422524159978\\
};

\addplot [color=mycolor2, draw=none, mark=square, mark options={solid, mycolor2}, forget plot]
 plot [error bars/.cd, y dir = both, y explicit]
 table[row sep=crcr, y error plus index=2, y error minus index=3]{%
3	29.4528179890485	1.86161475042722	1.86161475042722\\
1	15.7658242670572	1.02431120865221	1.02431120865221\\
4	20.8002356272759	1.19536716160624	1.19536716160624\\
2	13.3881418940295	0.821426018828489	0.821426018828489\\
9	13.213125331716	0.374624176657314	0.374624176657314\\
8	10.22076146703	0.284537910536524	0.284537910536524\\
7	14.4218541235825	0.402604898059883	0.402604898059883\\
5	7.99290085903223	0.460061471580098	0.460061471580098\\
10	13.2183031481995	0.367405040516955	0.367405040516955\\
6	8.24481585199791	0.460190132752167	0.460190132752167\\
17	8.4952737785251	0.216133413303992	0.216133413303992\\
15	7.69746547310033	0.205972196169537	0.205972196169537\\
13	13.3073050432551	0.384758996524126	0.384758996524126\\
11	22.2635004075045	0.680676221777662	0.680676221777662\\
18	8.22062578409	0.227873093254797	0.227873093254797\\
16	9.46688251368899	0.268045871379412	0.268045871379412\\
14	8.00510395782808	0.225132272439075	0.225132272439075\\
12	22.3335234533367	0.658168168905797	0.658168168905797\\
28	30.5890682087857	0.600526899667971	0.600526899667971\\
27	8.11158862530831	0.218428583352323	0.218428583352323\\
23	6.92060165069681	0.191752229425695	0.191752229425695\\
24	6.19459008559966	0.171760070078942	0.171760070078942\\
21	8.30357276533469	0.240956394385884	0.240956394385884\\
22	12.0956599766172	0.34251227129501	0.34251227129501\\
19	13.6373061440284	0.374087554719145	0.374087554719145\\
29	7.75476732940193	0.213386527168486	0.213386527168486\\
25	8.14552068300748	0.213088948368228	0.213088948368228\\
26	9.69613047747411	0.284877548048205	0.284877548048205\\
20	13.5449114796437	0.363750767837763	0.363750767837763\\
42	20.1919076598978	0.372360829752765	0.372360829752765\\
38	6.61164660903829	0.192152300939522	0.192152300939522\\
39	7.31005430784158	0.201175457624531	0.201175457624531\\
35	8.21992381851748	0.218260051506223	0.218260051506223\\
34	6.16900653238747	0.177710299194089	0.177710299194089\\
32	7.70716196562403	0.213754398961966	0.213754398961966\\
30	6.77547546316215	0.177704479007469	0.177704479007469\\
45	22.2075924088145	0.452071985997697	0.452071985997697\\
40	6.4658805993432	0.182282100243511	0.182282100243511\\
43	7.00957693868905	0.191881372011018	0.191881372011018\\
33	7.89270042542708	0.222810896957537	0.222810896957537\\
36	5.84918067519298	0.166356344162835	0.166356344162835\\
37	9.05265957840168	0.253295032659485	0.253295032659485\\
31	9.39391307783674	0.260175270787343	0.260175270787343\\
};

\end{axis}

\begin{axis}[%
width=5.75in,
height=1.35in,
at={(0in,0in)},
scale only axis,
xmin=0,
xmax=1,
ymin=0,
ymax=1,
axis line style={draw=none},
ticks=none,
axis x line*=bottom,
axis y line*=left,
legend style={legend cell align=left, align=left, draw=white!15!black}
]
\end{axis}
\end{tikzpicture}%
        \caption{Exact FGIs with coefficients from \cite{omelyan2003symplectic} (red circles) vs. Hessian-free FGIs with coefficients from Tab. \ref{tab:aFGIs} (blue squares).}
        \label{fig:SW_OMF_vs_HfFGI}
    \end{subfigure}

\medskip
    \begin{subfigure}[b]{\textwidth}
        \centering
        % This file was created by matlab2tikz.
%
%The latest updates can be retrieved from
%  http://www.mathworks.com/matlabcentral/fileexchange/22022-matlab2tikz-matlab2tikz
%where you can also make suggestions and rate matlab2tikz.
%
\definecolor{mycolor1}{rgb}{0.68627,0.22353,0.00000}%
\definecolor{mycolor2}{rgb}{0.00000,0.46275,0.68627}%
\begin{tikzpicture}

\begin{axis}[%
width=5in,
height=0.85in,
at={(0.5in,0.425in)},
scale only axis,
xmin=0,
xmax=46,
ymin=0,
ymax=40,
xlabel style={font=\color{white!15!black}},
xlabel={integrator ID in \cite{omelyan2003symplectic}}, 
ylabel style={font=\color{white!15!black}},
ylabel={$n_f/h(P_{\mathrm{acc}} = 90\%)$},
axis background/.style={fill=white},
axis x line*=bottom,
axis y line*=left,
legend style={legend cell align=left, align=left, draw=white!15!black}
]
\addplot [color=mycolor1, draw=none, mark=o, mark options={solid, mycolor1}, forget plot]
 plot [error bars/.cd, y dir = both, y explicit]
 table[row sep=crcr, y error plus index=2, y error minus index=3]{%
1	15.2696586332788	0.919060950372444	0.919060950372444\\
3	29.2288841364122	1.85284783693249	1.85284783693249\\
2	13.6515206586834	0.828322783054113	0.828322783054113\\
4	21.1148068161088	1.22789109141493	1.22789109141493\\
5	9.05184613207423	0.506926865027442	0.506926865027442\\
7	13.888554578916	0.384248663829379	0.384248663829379\\
8	10.3784476301091	0.27424402187989	0.27424402187989\\
9	13.2196563160468	0.360559765220287	0.360559765220287\\
6	7.95633743951507	0.426369803882757	0.426369803882757\\
10	12.9079513519028	0.363014807816982	0.363014807816982\\
11	21.8804134252911	0.654372290890404	0.654372290890404\\
13	19.3219142735609	0.52362441528779	0.52362441528779\\
15	7.89463734714292	0.213037648575953	0.213037648575953\\
17	8.9944857498538	0.231964984478359	0.231964984478359\\
12	21.8170134911382	0.633709786371046	0.633709786371046\\
14	8.30602886712715	0.241563710767613	0.241563710767613\\
16	9.85657087848872	0.286048413772114	0.286048413772114\\
18	8.18352496626363	0.214111964744243	0.214111964744243\\
19	13.3558699546649	0.385240206161552	0.385240206161552\\
21	9.11215456014357	0.26458918259975	0.26458918259975\\
22	21.4565229165122	0.585958565451891	0.585958565451891\\
23	7.06401738877355	0.198258432359226	0.198258432359226\\
24	6.28910357983195	0.171620165280875	0.171620165280875\\
27	9.71846485153024	0.257069706437277	0.257069706437277\\
28	30.3314110824349	0.570418345529013	0.570418345529013\\
20	13.5562352643627	0.389751817139045	0.389751817139045\\
25	8.05489851843086	0.210984535714461	0.210984535714461\\
26	9.91666042574882	0.270048924308903	0.270048924308903\\
29	7.54937281982576	0.214271813901425	0.214271813901425\\
30	6.64854780502786	0.181160262104329	0.181160262104329\\
32	8.0780322347681	0.22916743886293	0.22916743886293\\
34	6.10872656542301	0.171151913083069	0.171151913083069\\
35	9.55742233699655	0.250727654157864	0.250727654157864\\
38	19.9095174254973	0.535061427684787	0.535061427684787\\
39	16.3787022181084	0.44725586051216	0.44725586051216\\
42	11.949218958326	0.33536814374864	0.33536814374864\\
44	15.3563296910721	0.41816584390671	0.41816584390671\\
31	9.41024458002129	0.266649901934594	0.266649901934594\\
33	7.68301106239726	0.228636968864615	0.228636968864615\\
36	6.65600495672186	0.189850004579107	0.189850004579107\\
37	8.73908006235294	0.252287728071743	0.252287728071743\\
40	12.6503952363101	0.344232599199762	0.344232599199762\\
41	62.99330479035	1.86305366397729	1.86305366397729\\
43	16.3584145965364	0.426636535638442	0.426636535638442\\
45	22.8336364144654	0.640510201346684	0.640510201346684\\
};

\addplot [color=mycolor2, draw=none, mark=square, mark options={solid, mycolor2}, forget plot]
 plot [error bars/.cd, y dir = both, y explicit]
 table[row sep=crcr, y error plus index=2, y error minus index=3]{%
3	29.4528179890485	1.86161475042722	1.86161475042722\\
1	15.7658242670572	1.02431120865221	1.02431120865221\\
4	20.8002356272759	1.19536716160624	1.19536716160624\\
2	13.3881418940295	0.821426018828489	0.821426018828489\\
9	13.213125331716	0.374624176657314	0.374624176657314\\
8	10.22076146703	0.284537910536524	0.284537910536524\\
7	14.4218541235825	0.402604898059883	0.402604898059883\\
5	7.99290085903223	0.460061471580098	0.460061471580098\\
10	13.2183031481995	0.367405040516955	0.367405040516955\\
6	8.24481585199791	0.460190132752167	0.460190132752167\\
17	8.4952737785251	0.216133413303992	0.216133413303992\\
15	7.69746547310033	0.205972196169537	0.205972196169537\\
13	13.3073050432551	0.384758996524126	0.384758996524126\\
11	22.2635004075045	0.680676221777662	0.680676221777662\\
18	8.22062578409	0.227873093254797	0.227873093254797\\
16	9.46688251368899	0.268045871379412	0.268045871379412\\
14	8.00510395782808	0.225132272439075	0.225132272439075\\
12	22.3335234533367	0.658168168905797	0.658168168905797\\
28	30.5890682087857	0.600526899667971	0.600526899667971\\
27	8.11158862530831	0.218428583352323	0.218428583352323\\
23	6.92060165069681	0.191752229425695	0.191752229425695\\
24	6.19459008559966	0.171760070078942	0.171760070078942\\
21	8.30357276533469	0.240956394385884	0.240956394385884\\
22	12.0956599766172	0.34251227129501	0.34251227129501\\
19	13.6373061440284	0.374087554719145	0.374087554719145\\
29	7.75476732940193	0.213386527168486	0.213386527168486\\
25	8.14552068300748	0.213088948368228	0.213088948368228\\
26	9.69613047747411	0.284877548048205	0.284877548048205\\
20	13.5449114796437	0.363750767837763	0.363750767837763\\
42	20.1919076598978	0.372360829752765	0.372360829752765\\
38	6.61164660903829	0.192152300939522	0.192152300939522\\
39	7.31005430784158	0.201175457624531	0.201175457624531\\
35	8.21992381851748	0.218260051506223	0.218260051506223\\
34	6.16900653238747	0.177710299194089	0.177710299194089\\
32	7.70716196562403	0.213754398961966	0.213754398961966\\
30	6.77547546316215	0.177704479007469	0.177704479007469\\
45	22.2075924088145	0.452071985997697	0.452071985997697\\
40	6.4658805993432	0.182282100243511	0.182282100243511\\
43	7.00957693868905	0.191881372011018	0.191881372011018\\
33	7.89270042542708	0.222810896957537	0.222810896957537\\
36	5.84918067519298	0.166356344162835	0.166356344162835\\
37	9.05265957840168	0.253295032659485	0.253295032659485\\
31	9.39391307783674	0.260175270787343	0.260175270787343\\
};

\end{axis}

\begin{axis}[%
width=5.75in,
height=1.35in,
at={(0in,0in)},
scale only axis,
xmin=0,
xmax=1,
ymin=0,
ymax=1,
axis line style={draw=none},
ticks=none,
axis x line*=bottom,
axis y line*=left,
legend style={legend cell align=left, align=left, draw=white!15!black}
]
\end{axis}
\end{tikzpicture}%
        \caption{Hessian-free FGIs with coefficients from \cite{omelyan2003symplectic} (red circles) vs. Hessian-free FGIs with coefficients from Tab. \ref{tab:aFGIs} (blue squares). The coefficients from \cite{omelyan2003symplectic} with \mbox{ID = 41} yield $n_f/h(P_{\mathrm{acc}} = 90\%) \approx 63$ and are not displayed for visualization purposes.}
        \label{fig:SW_HfOMF_vs_HfFGI}
    \end{subfigure}
    
\medskip
    \begin{subfigure}[b]{\textwidth}
        \centering
        % This file was created by matlab2tikz.
%
%The latest updates can be retrieved from
%  http://www.mathworks.com/matlabcentral/fileexchange/22022-matlab2tikz-matlab2tikz
%where you can also make suggestions and rate matlab2tikz.
%
\definecolor{mycolor1}{rgb}{0.68627,0.22353,0.00000}%
\definecolor{mycolor2}{rgb}{0.00000,0.46275,0.68627}%
\begin{tikzpicture}

\begin{axis}[%
width=5in,
height=0.85in,
at={(0.5in,0.425in)},
scale only axis,
xmin=0,
xmax=46,
ymin=0,
ymax=40,
xlabel style={font=\color{white!15!black}},
xlabel={integrator ID in \cite{omelyan2003symplectic}}, 
ylabel style={font=\color{white!15!black}},
ylabel={$n_f/h(P_{\mathrm{acc}} = 90\%)$},
axis background/.style={fill=white},
axis x line*=bottom,
axis y line*=left,
legend style={legend cell align=left, align=left, draw=white!15!black}
]
\addplot [color=mycolor1, draw=none, mark=o, mark options={solid, mycolor1}, forget plot]
 plot [error bars/.cd, y dir = both, y explicit]
 table[row sep=crcr, y error plus index=2, y error minus index=3]{%
3	37.5575548457607	2.45457591744452	2.45457591744452\\
1	15.262460758336	0.929642536029298	0.929642536029298\\
4	27.2151299519401	1.5641894981709	1.5641894981709\\
2	13.7933007067276	0.836537173584071	0.836537173584071\\
9	16.4057694492425	0.446432785765842	0.446432785765842\\
8	12.4273021782154	0.330678153616308	0.330678153616308\\
7	16.2780617311959	0.464114229387948	0.464114229387948\\
5	8.14092157451435	0.47876029354041	0.47876029354041\\
10	16.9177180975559	0.45800072109165	0.45800072109165\\
6	8.83454932582519	0.440422855056221	0.440422855056221\\
17	9.61594199631758	0.262362824168096	0.262362824168096\\
15	9.45718888399195	0.263604287619511	0.263604287619511\\
13	14.2731906713846	0.395694101969564	0.395694101969564\\
11	22.3354032027909	0.642278277315072	0.642278277315072\\
18	8.79155343177971	0.257449921944628	0.257449921944628\\
16	11.2332732367659	0.317928692496452	0.317928692496452\\
14	9.6780616328294	0.261413906737932	0.261413906737932\\
12  22.2913410058656	0.744566801723968	0.744566801723968\\
28	36.8892040340136	0.743953577534106	0.743953577534106\\
27	8.1484747399349	0.236559900339058	0.236559900339058\\
23	6.27128956026646	0.17766033620257	0.17766033620257\\
24	6.59311708677394	0.180356699324769	0.180356699324769\\
21	10.3400839754977	0.279537823111878	0.279537823111878\\
22	13.1922272366448	0.365099587991039	0.365099587991039\\
19	13.9545053296672	0.382908870167723	0.382908870167723\\
29	8.10538314786555	0.224878365763968	0.224878365763968\\
25	9.81902618010487	0.260775906190468	0.260775906190468\\
26	11.5197433706884	0.322371848832388	0.322371848832388\\
20	13.386491197286	0.372359980556009	0.372359980556009\\
42	24.6288525084205	0.49770814857499	0.49770814857499\\
38	8.45048250970098	0.252908076552629	0.252908076552629\\
39	7.52815826179233	0.216405450208754	0.216405450208754\\
35	9.65653222204286	0.28600132770311	0.28600132770311\\
34	6.44988501236961	0.17609610941598	0.17609610941598\\
32	8.61944651801556	0.235104009431484	0.235104009431484\\
30	6.62301185509646	0.170594074120592	0.170594074120592\\
45	27.3845621930826	0.519020759891946	0.519020759891946\\
40	6.20415775303485	0.182784688597355	0.182784688597355\\
43	7.36457106894723	0.221449991867827	0.221449991867827\\
33	8.05899497135801	0.212977539001638	0.212977539001638\\
36	5.49625537907589	0.155825103355739	0.155825103355739\\
37	9.65293347338068	0.287294716273801	0.287294716273801\\
31	9.17620923816114	0.256818286229402	0.256818286229402\\
};

\addplot [color=mycolor2, draw=none, mark=square, mark options={solid, mycolor2}, forget plot]
 plot [error bars/.cd, y dir = both, y explicit]
 table[row sep=crcr, y error plus index=2, y error minus index=3]{%
3	29.4528179890485	1.86161475042722	1.86161475042722\\
1	15.7658242670572	1.02431120865221	1.02431120865221\\
4	20.8002356272759	1.19536716160624	1.19536716160624\\
2	13.3881418940295	0.821426018828489	0.821426018828489\\
9	13.213125331716	0.374624176657314	0.374624176657314\\
8	10.22076146703	0.284537910536524	0.284537910536524\\
7	14.4218541235825	0.402604898059883	0.402604898059883\\
5	7.99290085903223	0.460061471580098	0.460061471580098\\
10	13.2183031481995	0.367405040516955	0.367405040516955\\
6	8.24481585199791	0.460190132752167	0.460190132752167\\
17	8.4952737785251	0.216133413303992	0.216133413303992\\
15	7.69746547310033	0.205972196169537	0.205972196169537\\
13	13.3073050432551	0.384758996524126	0.384758996524126\\
11	22.2635004075045	0.680676221777662	0.680676221777662\\
18	8.22062578409	0.227873093254797	0.227873093254797\\
16	9.46688251368899	0.268045871379412	0.268045871379412\\
14	8.00510395782808	0.225132272439075	0.225132272439075\\
12	22.3335234533367	0.658168168905797	0.658168168905797\\
28	30.5890682087857	0.600526899667971	0.600526899667971\\
27	8.11158862530831	0.218428583352323	0.218428583352323\\
23	6.92060165069681	0.191752229425695	0.191752229425695\\
24	6.19459008559966	0.171760070078942	0.171760070078942\\
21	8.30357276533469	0.240956394385884	0.240956394385884\\
22	12.0956599766172	0.34251227129501	0.34251227129501\\
19	13.6373061440284	0.374087554719145	0.374087554719145\\
29	7.75476732940193	0.213386527168486	0.213386527168486\\
25	8.14552068300748	0.213088948368228	0.213088948368228\\
26	9.69613047747411	0.284877548048205	0.284877548048205\\
20	13.5449114796437	0.363750767837763	0.363750767837763\\
42	20.1919076598978	0.372360829752765	0.372360829752765\\
38	6.61164660903829	0.192152300939522	0.192152300939522\\
39	7.31005430784158	0.201175457624531	0.201175457624531\\
35	8.21992381851748	0.218260051506223	0.218260051506223\\
34	6.16900653238747	0.177710299194089	0.177710299194089\\
32	7.70716196562403	0.213754398961966	0.213754398961966\\
30	6.77547546316215	0.177704479007469	0.177704479007469\\
45	22.2075924088145	0.452071985997697	0.452071985997697\\
40	6.4658805993432	0.182282100243511	0.182282100243511\\
43	7.00957693868905	0.191881372011018	0.191881372011018\\
33	7.89270042542708	0.222810896957537	0.222810896957537\\
36	5.84918067519298	0.166356344162835	0.166356344162835\\
37	9.05265957840168	0.253295032659485	0.253295032659485\\
31	9.39391307783674	0.260175270787343	0.260175270787343\\
};

\end{axis}

\begin{axis}[%
width=5.75in,
height=1.35in,
at={(0in,0in)},
scale only axis,
xmin=0,
xmax=1,
ymin=0,
ymax=1,
axis line style={draw=none},
ticks=none,
axis x line*=bottom,
axis y line*=left,
legend style={legend cell align=left, align=left, draw=white!15!black}
]
\end{axis}
\end{tikzpicture}%
        \caption{Exact FGIs with coefficients from Tab. \ref{tab:aFGIs} (red circles) vs. Hessian-free FGIs with coefficients from Tab. \ref{tab:aFGIs} (blue squares).}
        \label{fig:HFGI_vs_HfFGI}
    \end{subfigure}
    \caption{Schwinger model. Estimation of the required number of force steps per unit trajectory $n_f/h$ to reach $P_\mathrm{acc} = 90\%$. Evaluations of the FG-term are counted as two force evaluations. The subfigures show a pair-wise comparison to Hessian-free FGIs with coefficients from Tab. \ref{tab:aFGIs} (blue squares) with (a) exact FGIs with coefficients from \cite{omelyan2003symplectic}, (b) Hessian-free FGIs with coefficients from \cite{omelyan2003symplectic}, (c) exact FGIs with coefficients from Tab. \ref{tab:aFGIs} (red circles).}
    \label{fig:SW_results}
\end{figure}
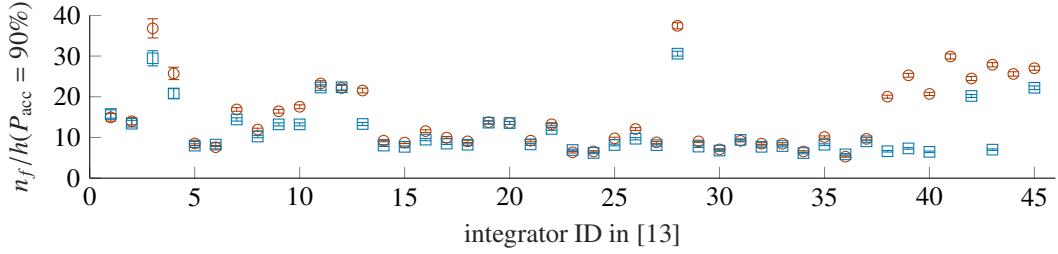
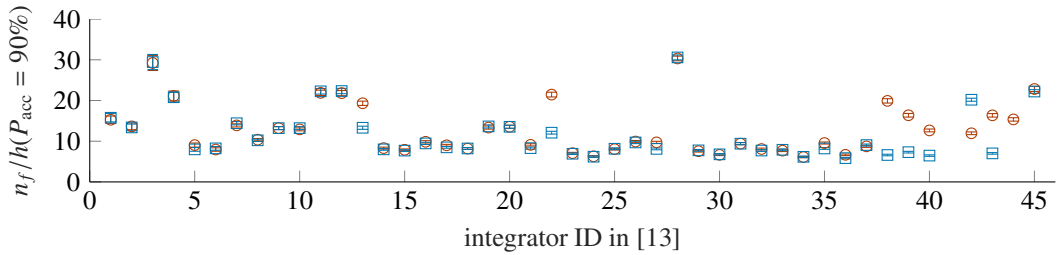
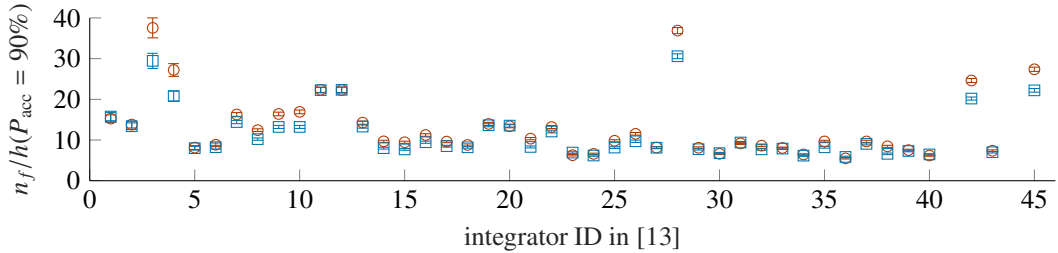

The two-dimensional Schwinger model is defined by the Hamiltonian 
\begin{align*}%\label{eq:Hamiltonian_2DSchwinger}
    \Ham(\p,\q) = \frac{1}{2} \sum\limits_{n=1,\mu=1}^{V,2} p_{n,\mu}^2 + \mathcal{S}_G([\bfQ]) + \mathcal{S}_F([\bfQ])
\end{align*}
with $V = L \times T$ the volume of the lattice and $\mathcal{S}_G([\bfQ])$ the pure gauge action and $\mathcal{S}_F([\bfQ])$ the fermion action. The links $\bfQ$ on the lattice are elements of the Abelian Lie group $\mathrm{U}(1)$. Consequently the links $Q_\mu(n)$ connecting the sites $n$ and $n+ \hat{\mu}$ and can be written in the form $Q_\mu(n) = \exp(i q_\mu(n)) \in \mathrm{U}(1)$, $q_{\mu}(n) \in [-\pi,\pi],\; \mu \in \{ x,t\}$ the respective space and time directions. The momenta $p_{n,\mu}$ are elements of the corresponding Lie algebra $\mathbb{R}$. 
Details on the action $\mathcal{S}$ can be found here \ref{app:Schwinger}.
For this Hamiltonian system, the force-gradient term $\hat{\C}$ has been derived in \cite{shcherbakov2017adapted}, allowing for a comparison of the exact FGIs in \cite{omelyan2003symplectic} and our Hessian-free ones. 
On a lattice of size $16\times16$ at pure gauge coupling $\beta=1$ and at bare mass parameter $m_0=0.352443$ \cite{Christian:2005yp},
we computed $2000$ trajectories of length $\tau = 1$ with $N=20,25,30,35,40$ number of steps using any set of coefficients from Tab.\ \ref{tab:aFGIs} and \cite[Tab.\ 2]{omelyan2003symplectic}.
Note that we chose for the simulations a large bare quark mass corresponding to $z=0.8$ in \cite[Tab.\ 2]{Christian:2005yp} to guarantee a stable integration for our setup without nested integration and without Hasenbusch mass-preconditioning in a large volume.
For all variants, we consider both the exact FG-term \eqref{eq:force-gradient_term} and the approximation \eqref{eq:Yin_derivation} resulting in Hessian-free FGIs. 
By measuring the dependence of the variance $\sigma^2$ of $\Delta H$ on the number of steps $N$ and using the relation to the acceptance rate \eqref{eq:acceptance_rate}, we estimated for each integrator the number of required force steps to reach $90\%$ acceptance rate.
In Fig.\ \ref{fig:SW_results}, a pair-wise comparison to our proposed integrators (Hessian-free FGIs with coefficients from Tab.\ \ref{tab:aFGIs}) is depicted. 
Fig.\ \ref{fig:SW_OMF_vs_HfFGI} shows that, for almost all variants, Hessian-free FGIs allow for a more efficient computational process than the FGIs from \cite{omelyan2003symplectic}. A major difference can be observed for those variants who are affected by an order reduction, e.g.~\eqref{eq:ADABADABADA}. This can be explained by the fact that in the relevant region of $P_{\mathrm{acc}} \approx 90\%$, second and fourth order integrators turn out to be more efficient for this lattice size. For increasing lattice sizes (see e.g.\ Section \ref{sec:QCD}), higher-order schemes become more efficient.
For simplicity, exact FGIs are often implemented using the approximation \eqref{eq:Yin_derivation}, see e.g.\ \cite{osborn2024tuning}. The coefficients derived in \cite{omelyan2003symplectic} do not incorporate the new error terms \eqref{eq:new_error_terms} that occur due to the approximation. Fig.\ \ref{fig:SW_HfOMF_vs_HfFGI} highlights that Hessian-free FGIs with coefficients from Tab.\ \ref{tab:aFGIs} usually perform better than Hessian-free FGIs with coefficients from \cite{omelyan2003symplectic}. 
In contrast, Fig.\ \ref{fig:HFGI_vs_HfFGI} underlines that the coefficients derived in this paper are only recommended for Hessian-free FGIs. Using the coefficients for exact FGIs leads to less efficient schemes, as expected. In this case, the coefficients from \cite{omelyan2003symplectic} should be used since they are optimized for exact FGIs.

\subsubsection{Four-dimensional lattice QCD simulations with Wilson fermions}\label{sec:QCD}
\begin{figure}[H]
    \centering
    % This file was created by matlab2tikz.
%
%The latest updates can be retrieved from
%  http://www.mathworks.com/matlabcentral/fileexchange/22022-matlab2tikz-matlab2tikz
%where you can also make suggestions and rate matlab2tikz.
%
\definecolor{mycolor1}{rgb}{0.68627,0.22353,0.00000}%
\definecolor{mycolor2}{rgb}{0.00000,0.46275,0.68627}%
\begin{tikzpicture}

\begin{axis}[%
width=2.8in,
height=1.57in,
at={(0.575in,0.45in)},
scale only axis,
xmode=log,
xmin=30,
xmax=80,
xtick={35,40,45,50,60,70},
xticklabels={$35$,$40$,$45$,$50$,$60$,$70$},
xminorticks=true,
xlabel style={font=\color{white!15!black}},
xlabel={$n_f \cdot N$},
ymode=log,
ymin=0.005,
ymax=1,
yminorticks=true,
ylabel style={font=\color{white!15!black}},
ylabel={$\sigma^2(\Delta H)$},
axis background/.style={fill=white},
legend style={at={(1.05,0.825)}, anchor=north west, legend cell align=left, align=left, draw=white!15!black}
]

\addplot [color=mycolor2,mark=triangle*]
 plot [error bars/.cd, y dir = both, y explicit]
 table[row sep=crcr, y error plus index=2, y error minus index=3]{%
 24  29.134991678994083  3.794530077316668  3.794530077316668 \\
 27  6.932824229197490  0.911316445017714  0.911316445017714 \\
 30  2.541253795765532  0.336814641970982  0.336814641970982 \\
 33  1.044119684818787  0.132607086347019  0.132607086347019 \\
 36  0.425543731090976  0.051501035805807  0.051501035805807 \\
 39  0.219460552018306  0.034688295380241  0.034688295380241 \\
 42  0.107048385435466  0.016645792792137  0.016645792792137 \\ 
 45  0.065209177425111  0.008672929205028  0.008672929205028 \\ 
 48  0.036112351271265  0.004709946392832  0.004709946392832 \\ 
 51  0.025022034258571  0.003206031504063  0.003206031504063 \\ 
 54  0.014577080706361  0.002170578953453  0.002170578953453 \\ 
 57  0.007353427564932  0.001053434753774  0.001053434753774 \\ 
 60  0.005775131886391  0.000928916169344  0.000928916169344 \\ 
};
\addlegendentry{BADAB ($\mathrm{Eff} \approx 16.96$)}

\addplot [color=mycolor2,mark=pentagon*]
 plot [error bars/.cd, y dir = both, y explicit]
 table[row sep=crcr, y error plus index=2, y error minus index=3]{%
24  57.142071005917046  9.000787444431994  9.000787444431994  \\
28  5.482180288461535  0.688495758424982  0.688495758424982  \\
32  0.802287661150148  0.132267218017337  0.132267218017337  \\
36  0.130620106686391  0.019823032897364  0.019823032897364  \\
40  0.041682357850851  0.005115571695518  0.005115571695518  \\
48  0.007400132278763  0.001215839699153  0.001215839699153 \\  
52  0.003633145982803  0.000502756292364  0.000502756292364  \\  
56  0.002046006996117  0.000282742766348  0.000282742766348 \\  
60  0.000842628453217  0.000136549831478  0.000136549831478\\  
};
\addlegendentry{ABADABA ($\mathrm{Eff} \approx 26.19$)}

\addplot [color=mycolor2,mark=diamond*]
 plot [error bars/.cd, y dir = both, y explicit]
 table[row sep=crcr, y error plus index=2, y error minus index=3]{%
30  3.683748018398671  0.562789954518374  0.562789954518374 \\   
35  0.682886141034578  0.095936771137687  0.095936771137687 \\   
40  0.115237649973955  0.015558526205245  0.015558526205245 \\  
45  0.031565144311206  0.004136636755812  0.004136636755812 \\ 
50  0.007321411049371  0.000887849065550  0.000887849065550 \\ 
55  0.004062418176775  0.000605931019165  0.000605931019165 \\  
60  0.001884610514053  0.000223675280399  0.000223675280399 \\  
65  0.000927596715052  0.000126140538482  0.000126140538482 \\  
70  0.000515297010910  0.000062879867925  0.000062879867925 \\  
75  0.000291804075592  0.000038623683033  0.000038623683033 \\  
80  0.000125911934470  0.000017253011194  0.000017253011194 \\
};
\addlegendentry{BABADABAB ($\mathrm{Eff} \approx 24.57$)}

\addplot [color=mycolor1,mark=triangle*, mark options={solid, rotate=180, mycolor1}]
 plot [error bars/.cd, y dir = both, y explicit]
 table[row sep=crcr, y error plus index=2, y error minus index=3]{%
30  8.773805981878706  1.125620078447276  1.125620078447276 \\   
35  1.176030293546598  0.193849998952696  0.193849998952696 \\   
40  0.146297861907359  0.020778067050962  0.020778067050962 \\   
45  0.015072474984283  0.002101513015934  0.002101513015934 \\  
55  0.001126513034461  0.000136726137419  0.000136726137419 \\  
60  0.000544314772032  0.000069490214627  0.000069490214627 \\  
65  0.000219151757763  0.000029467877464  0.000029467877464 \\  
70  0.000102442989436  0.000013721724358  0.000013721724358 \\  
75  0.000068046722448  0.000009117465730  0.000009117465730 \\
80  0.000038834656059  0.000006046774706  0.000006046774706 \\
};
\addlegendentry{BABABABABAB ($\mathrm{Eff} \approx 59.26$)}

\addplot [color=mycolor1,mark=square*]
 plot [error bars/.cd, y dir = both, y explicit]
 table[row sep=crcr, y error plus index=2, y error minus index=3]{%
42  6.905804942809729  0.885941434228681  0.885941434228681 \\  
49  0.632935758783284  0.069521967576085  0.069521967576085 \\  
56  0.080226390324519  0.013206746044541  0.013206746044541 \\ 
63  0.009745777840237  0.001527961169086  0.001527961169086 \\  
70  0.001349226046366  0.000167882842307  0.000167882842307 \\  
77  0.000369107081130  0.000084139145755  0.000084139145755 \\  
84  0.000076776556797  0.000009759171880  0.000009759171880 \\  
91  0.000029388800156  0.000006965638167  0.000006965638167 \\  
98  0.000008484021662  0.000001129969680  0.000001129969680 \\
};
\addlegendentry{BABABABABABABAB ($\mathrm{Eff} \approx 1.40$)}

\addplot [color=black, dashed, forget plot]
  table[row sep=crcr]{%
1	0.0631631\\
70	0.0631631\\
};
\node[fill=white] at (axis cs: 67.5,0.0631631) {$\langle P_{\mathrm{acc}} \rangle = 90\%$};

\addplot [color=black, dashed, forget plot]
  table[row sep=crcr]{%
1	0.312076\\
70	0.312076\\
};
\node[fill=white] at (axis cs: 67.5,0.312076) {$\langle P_{\mathrm{acc}} \rangle = 78\%$};

\end{axis}

\begin{axis}[%
width=6.15in,
height=2.1in,
at={(0in,0in)},
scale only axis,
xmin=0,
xmax=1,
ymin=0,
ymax=1,
axis line style={draw=none},
ticks=none,
axis x line*=bottom,
axis y line*=left,
legend style={legend cell align=left, align=left, draw=white!15!black}
]
\end{axis}
\end{tikzpicture}%
    \caption{Lattice QCD. Variance of $\Delta H$ vs.\ number of force evaluations per trajectory $n_f \cdot N$ for a selected number of Hessian-free FGIs, namely the three best performing Hessian-free FGIs (blue lines), as well as the most efficient non-gradient scheme for order $p \in \{4,6\}$ (red lines). The most efficient non-gradient scheme of order $p=2$ is not included for visualization purposes as it achieves $\langle P_{\mathrm{acc}}\rangle_\mathrm{opt}$ at $n_f \cdot N \approx 140$.}
    \label{fig:lattice_QCD_results}
\end{figure}
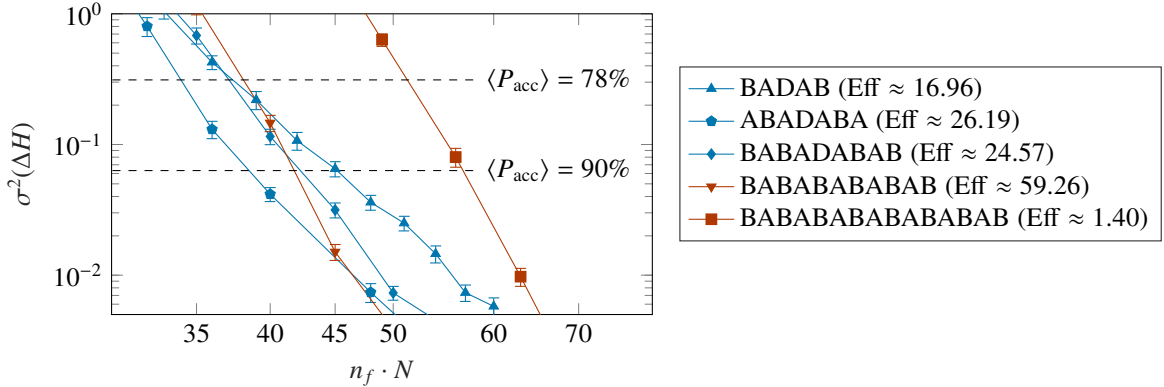

We consider gauge field simulations in lattice quantum chromodynamics (QCD) on a four-dimensional lattice of size $V = T \times L^3$ with lattice spacing $a$. Here, the Hamiltonian reads 
\begin{align}\label{eq:QCD_Hamiltonian}
    \Ham([\bfP],[\bfQ]) = \frac{1}{2}\sum\limits_{x,\mu} \mathrm{tr}(P_{x,\mu}^2) + \mathcal{S}_G([\bfQ]) + \mathcal{S}_{F}([\bfQ]),
\end{align}
where $\mathcal{S}_G([\bfQ])$ denotes the Wilson gauge action, $\mathcal{S}_{F}([\bfQ])$ denotes the fermion action (for details, see \ref{app:latticeQCD}).
Here, the links $Q_{x,\mu}$, connecting the sites $x$ and $x + a\hat{\mu}$, are elements of the non-Abelian matrix Lie group $\mathrm{SU}(3)$ of unitary matrices with determinant one.
The scaled momenta $iP_{x,\mu}$ are elements of the associated Lie algebra $\mathfrak{su}(3)$ of traceless and anti-Hermitian matrices and occur in the formulation of the equations of motion.
As an initial test in lattice QCD simulations, we consider an ensemble used in \cite{Knechtli2022} generated with two dynamical non-perturbatively
$\mathcal{O}(a)$ improved Wilson quarks at a mass equal to
half of the physical charm quark.
Starting from a thermalized configuration, we computed 100 trajectories for $\tau=2.0$ and varying step size $h = \tau/N$ for different values of $N$ on a $48 \times 24^3$ lattice with gauge coupling $\beta = 5.3$ and hopping parameter $\kappa = 0.1327$ for all Hessian-free FGIs discussed in this paper.
For the simulations, we have put all forces on a single time-scale of integration.\footnote{We also performed simulations using nested integrators that employ a smaller step size to the gauge force. At these lattice parameters, however, we were not able to observe any significant improvement in the acceptance probability by applying nested integration techniques, see e.g. \cite{sexton1992hamiltonian,URBACH200687}.}
The results of selected integrators are depicted in Fig.\ \ref{fig:lattice_QCD_results}. 

The results emphasize that the use of fourth-order integrators and in particular the Hessian-free FGI \eqref{eq:ABADABA} results in the most efficient computational process with $n_f \cdot N \approx 34$ to achieve the optimal acceptance rate $\langle P_\mathrm{acc}\rangle_\mathrm{opt} = 78\%$. In contrast, the most efficient non-gradient scheme \eqref{eq:BABABABABAB} demands $n_f \cdot N \approx 38$. Moreover, the Hessian-free FGI \eqref{eq:BADAB} is a very efficient choice ($n_f \cdot N \approx 37$). The integrators \eqref{eq:ABADABA} ($\mathrm{Eff} \approx 26.19$) and \eqref{eq:BADAB} ($\mathrm{Eff} \approx 16.96$) are not among the most efficient integrators according to the efficiency measure \eqref{eq:efficiency_measure}. In lattice QCD, however, we are only interested in the region of $0.01 \leq \sigma^2 \leq 1$. In this region, the limiting factor in the step size $h$ usually comes from instabilities \cite{instabilityMD} so that the performance of the integrators strongly depends on their stability properties. For the non-gradient schemes, a linear stability analysis is already available \cite{BlanesStability}. First results on the linear stability of FGIs support the claim that Hessian-free FGIs with a larger stability domain are beneficial for application in lattice QCD simulations. Moreover, we emphasize that Hessian-free FGIs are not symplectic and thus will in general have a linear energy drift of size $\mathcal{O}(\tau h^{\max\{4,p\}})$ while exact FGIs are symplectic, resulting in a bounded energy error without any drift. Since the trajectory length $\tau$ is rather short, the energy drift is negligible. The linear stability and energy conservation of Hessian-free FGIs is subject of a following paper in progress.

For the best-performing Hessian-free FGI \eqref{eq:ABADABA} and the most efficient non-gradient integrator \eqref{eq:BABABABABAB}, we tuned the step size $h$ to achieve an acceptance rate of $P_{\mathrm{acc}} \geq 90\%$, resulting in the step size $h=0.2$ for both integrators, and computed 2000 trajectories. The results of both setups are summarized in Table \ref{tab:QCDResults}, including the cost measure 
$$ \mathrm{cost} = \frac{n_f \cdot N}{P_{\mathrm{acc}} \cdot \tau},$$
as well as estimates of the integrated autocorrelation times of the topological charge $Q_0$ measured at Wilson flow time $t_0$, and the Wilson flow reference scale $t_0/a^2$.
For both quantities we observe similar autocorrelation times, despite the lower acceptance rate of the Hessian-free force-gradient integrator.
By investigating the cost measure, it turns out that the Hessian-free FGI only demands approximately $84\%$ of the computational cost compared to the non-gradient integrator. 
Despite the simplicity of the ensemble, the results show the potential of Hessian-free FGIs for lattice QCD simulations. In more complicated ensembles, the performance usually benefits from applying nested integration techniques \cite{sexton1992hamiltonian,URBACH200687} to employ different step sizes for different parts of the forces. The investigation of Hessian-free FGIs inside nested integrators is subject to future research focusing on multiple time stepping techniques.
\begin{table}[H]
    \centering
    \begin{tabular}{l r r r r r}
        \toprule
        ID & $n_f \cdot N$ & $P_\mathrm{acc}$ & $\tau_{\mathrm{int}}(t_0)$ [MDU] & $\tau_{\mathrm{int}}(Q_0)$ [MDU] & cost \\\midrule
        BABABABABAB & $50$ & $97.5\%$ & $37.40(12.66)$ & $22.91(6.52)$ & $25.64$ \\
        ABADABA & $40$ & $92.3\%$ & $28.15(9.24)$ & $22.05(6.31)$ & $21.67$ \\
        \bottomrule
    \end{tabular}
    \caption{Lattice QCD. Comparison of the performance of two tuned setups based on 2000 trajectories of length $\tau = 2.0$. Both setups use a step size of $h=0.2$. The results contain the number of force evaluations per trajectory, the acceptance rate $P_{\mathrm{acc}}$, integrated autocorrelation times $\tau_{\mathrm{int}}(t_0)$ and $\tau_{\mathrm{int}}(Q_0)$ in molecular dynamics units (MDUs), and a measure for the computational cost.}
    \label{tab:QCDResults}
\end{table}

\section{Conclusion and outlook}
In this work, we generalized the approximation proposed in \cite{wisdom1996symplectic,yin2011improving}, resulting in a new class of Hessian-free FGIs. The new integrators are applicable to separable Hamiltonian systems of the form \eqref{eq:Hamiltonian_general}, as well as to general second order ODEs $\ddot{y} = f(y)$. 
In contrast to exact FGIs~\cite{omelyan2003symplectic}, the Hessian-free variants do not need for the Hessian of the potential energy $\V$ and replace it by another force evaluation. At the price of additional error terms, that have been derived explicitly in this work, one is able to save costly evaluations of the force-gradient term, e.g.\ in molecular dynamics simulations. 
Hessian-free FGIs are time-reversible, volume-preserving and satisfy the closure property, but they are no longer symplectic. Consequently, they do not preserve a nearby shadow Hamiltonian. In general, this results in a linear energy drift of size $\mathcal{O}(th^p)$.
A complete classification of Hessian-free FGIs up to $s=11$ stages shows that some variants undergo only small changes in the time coefficients compared to the respective variants in \cite{omelyan2003symplectic}, whereas the time coefficients of other variants change significantly, e.g.\ if an order reduction appears due to the additional error terms. 
Numerical tests emphasize the computational efficiency of the proposed numerical integration schemes for different applications, making them a valuable class of integration schemes for a wide range of applications.
The efficiency of higher order integration schemes strongly depends on the system size so that Hessian-free FGIs of order $p>4$ may become efficient for simulating lattice quantum field theories of larger volumes \cite{takaishi2006}.
Hence a next natural step is to investigate Hessian-free FGIs with $s>11$ stages, aiming for more efficient sixth-order integrators.
In the context of lattice QCD, the numerical results indicate that the stability properties of the integrators are crucial for an efficient computational process. Thus we will perform a linear stability analysis of Hessian-free and exact FGIs to validate which variants are superior in terms of stability. 
Furthermore, a refined analysis on the energy conservation of Hessian-free FGIs in lattice QCD will contribute in designing more efficient numerical integration schemes for simulating quantum field theories on the lattice.
The efficiency measure used in this work assumes that the error terms are equally dominant. However, this simplifying assumption does not hold in general. Minimizing a problem-specific weighted norm results in even more efficient integrators \cite{takaishi2006}. Therefore, we will develop a tuning model that numerically measures the error terms to obtain appropriate weights for the weighted norm. Since we derived the additional error terms for Hessian-free FGIs in this work, this enables tuning procedures similar to \cite{takaishi2006} and may enhance other existing tuning models, e.g.\ \cite{osborn2024tuning}.
In many applications, separable Hamiltonian systems with a potential splitting $\Ham(\p,\q) = \T(\p) + \sum\nolimits_{m=1}^N \V^{\{m\}}(\q)$ occur. Here, the potential part is split with respect to, for example, stiffness, nonlinearity, dynamical behavior, and evaluation cost. In lattice QCD, for example, the fermion part of the action is characterized by slow dynamics and expensive evaluation cost, whereas the gauge part is fast and cheap to evaluate. This multirate potential can be exploited by means of nested integration, see e.g.\ \cite{shcherbakov2017adapted,sexton1992hamiltonian,URBACH200687}. We aim at designing a general framework of nested integration schemes and providing an order theory including explicit formulae of the error terms. 

\section*{CRediT authorship contribution statement}
\textbf{Kevin Schäfers:} Conceptualization, Methodology, Software, Validation, Investigation, Writing - Original Draft,\ 
\textbf{Jacob Finkenrath:} Software, Validation, Investigation, Writing - Original Draft,\ 
\textbf{Michael Günther:} Writing - Review \& Editing, Supervision, Funding acquisition,\ 
\textbf{Francesco Knechtli:} Writing - Review \& Editing, Supervision, Project administration, Funding acquisition.

\section*{Declaration of competing interest}
The authors declare that they have no known competing financial interests or personal relationships that could have appeared to influence the work reported in this paper.

\section*{Data availability}
An extended version of the openQCD package \cite{openQCD} including the proposed class of integrators is opensource and available online at \url{https://github.com/KevinSchaefers/openQCD_force-gradient}.

\section*{Acknowledgements}
This work is supported by the German Research Foundation (DFG) research unit FOR5269 "Future methods for studying confined gluons in QCD".\ 
The lattice QCD simulations were carried out on the PLEIADES cluster at the University of Wuppertal, which was supported by the DFG and the Federal Ministry of Education and Research (BMBF), Germany.\ 
We thank Stefan Schaefer for advice on the openQCD code, as well as Mike Peardon and Tomasz Korzec for discussions.
We gratefully acknowledge Robert I. McLachlan for his helpful insights on an issue in the initial version of the manuscript.

%% The Appendices part is started with the command \appendix;
%% appendix sections are then done as normal sections
\appendix

\section{Details on the action for the two-dimensional Schwinger model}\label{app:Schwinger}
The action $\mathcal{S}$ is given by the sum of the gauge part 
\begin{align*}%\label{eq:Schwinger_gauge-action}
    \mathcal{S}_G([\bfQ]) = \beta \sum_{n=1}^V (1 - \mathfrak{Re}(\mathbb{P}(n)))
\end{align*}
summing up over all plaquettes
$$\mathbb{P}(n) = Q_1(n)Q_2(n+\hat{1})Q_1^\dagger (n+\hat{2}) Q_2^\dagger(n),$$
scaled by a coupling constant $\beta$, and the fermion part 
\begin{equation*}%\label{eq:Schwinger_fermion-action}
    \mathcal{S}_F([\bfQ]) = \eta^\dagger \left(D^\dagger D\right)^{-1} \eta,
\end{equation*}
with pseudofermion field $\eta$ and Wilson-Dirac operator 
$$ D_{n,m} = (2+m_0) \delta_{n,m} - \frac{1}{2} \sum_{\mu=1}^2 \left((1-\sigma_\mu) Q_\mu(n) \delta_{n,m-\hat{\mu}} + (1+ \sigma_\mu) Q_\mu^\dagger(n-\hat{\mu}) \delta_{n,m+\hat{\mu}} \right),$$
where $m_0$ denotes the mass parameter, $\delta_{n,m}$ the Kronecker delta acting on the pseudofermion field via 
$$\sum_{m=1}^V \delta_{n,m} \eta(m) = \eta(n),$$ 
and $\sigma_\mu$ the Pauli matrices 
$$\sigma_1 = \begin{pmatrix}
    0 & 1 \\ 1 & 0
\end{pmatrix}, \quad \sigma_2 = \begin{pmatrix}
    0 & -i \\ i & 0
\end{pmatrix},$$
which build a basis of the Lie group $\mathrm{U}(1)$. 

\section{Details on the Hamiltonian for lattice QCD}\label{app:latticeQCD}
Any matrix representation $P_{x,\mu}$ can be written as $P_{x,\mu} = p^i \bfT_i$ where the generators $\bfT_i,\, i=1,\ldots,8$, are given by the Gell-Mann matrices \cite{georgi2000lie}. With $\p_{x,\mu} \coloneqq (p_1,\ldots,p_8)^\top$, a straight-forward computation shows that 
$$\mathrm{tr}(P_{x,\mu}^2) = \langle \p_{x,\mu},\p_{x,\mu}\rangle,$$
holds. Consequently, the kinetic energy in \eqref{eq:QCD_Hamiltonian} can be written as 
$$ \T(\p) = \frac{1}{2} \sum\limits_{x,\mu} \langle \p_{x,\mu},\p_{x,\mu} \rangle$$
which is of the form \eqref{eq:T} with $\bfM = \mathrm{Id}$.
The action $\mathcal{S} = \mathcal{S}_G + \mathcal{S}_F$ is defined on a four dimensional lattice with volume $V$ and consists of the gauge part 
\begin{equation*}
    \mathcal{S}_G([\bfQ]) = \sum_x \sum_{\mu < \nu} \beta \left( 1 - \frac{1}{3} \mathfrak{Re}\left( \mathrm{tr}\left( \mathbb{P}_{x,\mu\nu}([\bfQ])\right)\right)\right)
\end{equation*}
with gauge coupling $\beta$ and plaquettes
$$ \mathbb{P}_{x,\mu\nu}([\bfQ]) = Q_{x,\mu} Q_{x + a\hat{\mu}} Q_{x+a\hat{\nu},\mu}^\dag Q_{x,\nu}^\dag,$$
and the fermion part
\begin{equation*}
    \mathcal{S}_F([\bfQ]) = \langle \phi_1, (D_{eo}D_{eo}^\dag+\mu^2)^{-1} \phi_1 \rangle + \langle \phi_2, (1 + \mu^2 (D_{eo}D_{eo}^\dag)^{-1} )\phi_2 \rangle + \sum_x \textrm{tr}\,\textrm{log} D^2_{x,x} 
\end{equation*}
with complex-valued pseudofermion fields $\phi_1$ and $\phi_2$.
The fermion part is decomposed using even-odd reduction in combination
with one Hasenbusch mass preconditioning term \cite{Hasenbusch:2001ne} with shift parameter $\mu$,
such that the determinant is given by
\begin{equation*}
\det D^2 = \frac{\det (D_{eo} D_{eo}^\dagger+ \mu^2 \cdot \mathrm{Id} )} { \det (\mathrm{Id} + \mu^2 (D_{eo} D_{eo}^\dagger )^{-1})} \cdot \prod_x \det D^2_{x,x}
\end{equation*}
with the even-odd reduced Dirac operator $D_{eo} = \mathrm{Id}-D_{e,e}^{-1} D_{e,o}D_{o,o}^{-1}D_{o,e}$ and the block diagonal terms $D_{x,x}$. The clover improved Wilson Dirac operator can be represented as a complex matrix of size $12 \times V$ with next neighbour interactions.
For a more detailed discussion, we refer to \cite{Frommer:2013fsa}.

\bibliographystyle{elsarticle-num} 
\bibliography{refs}
\end{document}